\documentstyle[11pt]{article}

\catcode`@=11
\def\UU{\leavevmode\setbox0\hbox{h}\dimen@\ht0\advance\dimen@-1ex%
  \rlap{\raise.67\dimen@\hbox{\char'27}}U}
\def\TT{\leavevmode\setbox0\hbox{h}\dimen@\ht0\advance\dimen@-1ex%
  \rlap{\raise.67\dimen@\hbox{\char'27}}T}
\def\lTT{\leavevmode\setbox0\hbox{h}\dimen@\ht0\advance\dimen@-1ex%
  \rlap{\lower3.1\dimen@\hbox{\char'27}}T}
\catcode`@=12

\newfont{\myfont}{cmbx10}

\title{ {\bf Algorithm xxx --- ORTHPOL: A package of routines
for generating orthogonal \hbox{polynomials and Gauss-type quadrature
rules}\thanks{Work supported, in part, by
National Research Foundation grants, most recently by grant DMS-9023403.}}}

\author{WALTER GAUTSCHI}

\date{\ }
\begin{document}

\maketitle

\vspace{-.75in}

\bigskip

\small
ABSTRACT.
A collection of subroutines and examples of their uses, as well as
the underlying numerical methods, are described for generating
orthogonal polynomials relative to arbitrary weight functions.
The object of these routines is to produce the coefficients in the
three-term recurrence relation satisfied by the orthogonal
polynomials.
Once these are known, additional data can be generated, such as
zeros of orthogonal polynomials and Gauss-type quadrature rules,
for which routines are also provided.
\normalsize
\bigskip


\noindent
CR {\it Classification Scheme}. G.1.2, G.1.4, G.4.
\bigskip

\noindent
{\small 1.  INTRODUCTION}
\bigskip

Classical orthogonal polynomials, such as those of Legendre, Chebyshev,
Laguerre and Hermite, have been used for purposes of approximation in
widely different disciplines and over a long period of time.
Their popularity is due in part to the ease with which they can be employed,
and in part to the wealth of analytic results known for them.
Widespread use of nonclassical orthogonal
polynomials, in contrast, has been impeded by a lack of
effective and generally applicable constructive methods.
The present set of computer routines has been developed over the past
ten years in the hope of remedying this impediment and of
encouraging the use of nonstandard orthogonal polynomials.
A number of applications indeed have already been made, for example
to numerical quadrature (Cauchy principal value integrals with
coth-kernel [38], Hilbert transform of Jacobi weight functions [37],
integration over half-infinite intervals [28], rational Gauss-type
quadrature [30, 31]), to moment-preserving
spline approximation [21,35,11], to the summation of slowly convergent
series [26,27], and, perhaps most notably, to the proof of the
Bieberbach conjecture [24].

In most applications, orthogonality is with respect to a positive weight
function, $w$, on a given interval or union of intervals, or with respect to
positive weights, $w_i$,
concentrated on a discrete set of points, $\{ x_i \}$, or a
combination of both.
For convenience of notation, we subsume all these cases under the notion
of a positive measure, $d \lambda$, on the real line
{\myfont R}.
That is, the respective inner product is written as a
Riemann-Stieltjes integral,
$$
(u,v) = \int_{\mbox{\myfont R} } u(t) v(t) d \lambda (t),
\eqno(1.1)
$$
where the function $\lambda (t)$ is the indefinite integral of
$w$ for the continuous part, and a step function with jumps $w_i$ at $x_i$
for the discrete part.
We assume that (1.1) is meaningful whenever $u,v$ are polynomials.
There is then defined a unique set of (monic) orthogonal polynomials,
\renewcommand{\arraystretch}{2}
$$
\begin{array}{c}
{\displaystyle
\pi_k (t) = t^k + ~
\mbox{lower-degree~terms}, ~~~~ k = 0, 1, 2,\ldots, } \\
{\displaystyle
( \pi_k , \pi_{\ell} ) = 0 ~~ \mbox{if} ~~ k \neq \ell .}
\end{array}
\eqno(1.2)
$$
\renewcommand{\arraystretch}{1}We 
speak of ``continuous'' orthogonal
polynomials if the support of $d \lambda$ is an interval,
or a union of intervals, of ``discrete'' orthogonal polynomials if the
support of $d \lambda$ consists of a discrete set of points, and of
orthogonal polynomials of ``mixed type'' if
the support of $d \lambda$ has both a continuous and discrete part.
In the first and last case, there are infinitely many orthogonal polynomials ---
one for each degree ---, whereas in the second case there are exactly $N$
orthogonal polynomials,
$\pi_0 , \pi_1 ,\ldots, \pi_{N-1}$, where $N$ is the number of support
points.
In all cases, we denote the polynomials by
$\pi_k ( \cdot ) = \pi_k(~ \cdot ~; d \lambda )$, or
$\pi_k ( ~\cdot~; w )$,
if we want to indicate
their dependence on the measure $d \lambda$ or weight function $w$,
and use similar notations for other quantities depending on
$d \lambda$ or $w$.

It is a distinctive feature of orthogonal polynomials, compared
to other orthogonal systems, that they satisfy a three-term recurrence
relation,
\renewcommand{\arraystretch}{2}
$$
\begin{array}{c}
{\displaystyle
\pi_{k+1} (t) = (t - \alpha_k ) \pi_k (t) - \beta_k \pi_{k-1} (t) , ~~
k = 0, 1, 2 ,\ldots, }  \\
{\displaystyle
\pi_0 (t) = 1 , ~~ \pi_{-1} (t) = 0 ,}
\end{array}
\eqno(1.3)
$$
\renewcommand{\arraystretch}{1}with coefficients
$\alpha_k = \alpha_k ( d \lambda ) \in \mbox{\myfont R}$, $\beta_k = \beta_k ( d \lambda ) > 0$ that
are uniquely determined by the measure $d \lambda$.
By convention, the coefficient $\beta_0$, which multiplies
$\pi_{-1} = 0$ in (1.3), is defined by
$$
\beta_0 = \beta_0 ( d \lambda ) = \int_{\mbox{\myfont R} } d \lambda (t) .
\eqno(1.4)
$$
The knowledge of these coefficients is absolutely indispensable
for any sound computational use and application of orthogonal
polynomials [19,25].  One of the principal objectives of the present package
is precisely to provide routines for generating these coefficients.
Routines for
related quantities are also provided, such as
Gauss-type quadrature weights and nodes,
hence also zeros of orthogonal polynomials.

Occasionally (e.g., in [21,35,11,30,31]),
one needs to deal with indefinite (i.e., sign-changing) measures
$d \lambda$.
The positivity of the $\beta_k$ is then no longer guaranteed, indeed not
even the existence of all orthogonal polynomials.
Nevertheless, our methods can still be formally applied, albeit at the risk of
possible breakdowns or instabilities.

There are basically four methods used here to generate recursion
coefficients: (i) {\it Methods based on explicit formulae}.
These relate to classical orthogonal polynomials, and are implemented in
the routine {\tt recur} of \S2.  (ii) {\it Methods based on moment information}.
These are dealt with in \S3, and are represented by a single routine,
{\tt cheb}.  Its origin can be traced back to work of Chebyshev on discrete
least squares approximation.  (iii) {\it Bootstrap methods based on inner
product formulae
for the coefficients, and orthogonal reduction methods}.
We have attributed the idea for the former method to Stieltjes and referred to it
in [19] as the Stieltjes procedure.  The prototype is the routine {\tt sti}
in \S4, applicable for discrete orthogonal polynomials.
An alternative routine is {\tt lancz},
which accomplishes the same purpose, but uses the {\it method of
Lanczos}. Either of these routines can be used in {\tt mcdis},
which applies to continuous as well as mixed-type orthogonal polynomials.
In contrast to all previous
routines, {\tt mcdis} uses a discretization process, and thus furnishes only approximate answers
whose accuracies can be controlled by the user.
The routine, however, is by far the most sophisticated and flexible routine
in this package, one
that requires, or can greatly benefit from,
ingenuity of the user.
The same kind of discretization is also applicable to moment-related
methods, yielding the routine {\tt mccheb}.
(iv)  {\it Modification algorithms}.  These are routines
generating recursion coefficients for measures modified by a
rational factor, 
utilizing the recursion coefficients of the
original measure, which are assumed to be known.
They can be thought of as algorithmic implementations of the
Christoffel, or generalized Christoffel, theorem and are
incorporated in the routines
{\tt chri} and {\tt gchri} of \S5.
An important application of all these routines is made
in \S6, where routines are provided that generate the
weights and nodes of quadrature rules of Gauss, Gauss-Radau, and
Gauss-Lobatto type.

Each routine has a single-precision and double-precision version with
similar names, except for the prefix ``{\tt d}'' in double-precision
procedures.
The latter are generally a straightforward translation of the former.
An exception is the routine {\tt dlga} used in {\tt drecur} for computing
the logarithm of the gamma function, which employs a different method than
the single-precision companion routine {\tt alga}.

All routines of the package have been checked for ANSI conformance,
and have been tested on two computers: the
Cyber 205 and a Sun 4/670 MP workstation. The former has machine precisions
$\epsilon^{s} \approx 7.11 \times 10^{-15}$, $\epsilon^{d} \approx 5.05
\times 10^{-29}$ in single and double precision, respectively, while
the latter has $\epsilon^{s} \approx 5.96 \times 10^{-8}$, $\epsilon^{d}
\approx 1.11 \times 10^{-16}$. The Cyber 205 has a large floating-point
exponent range, extending from approximately $-8617$ to +8645 in single
as well as double precision, whereas the Sun 4/670 has the rather
limited exponent range [--38, 38] in single precision, but a larger
range, [--308, 308], in double precision. All output cited relates
to work on the Cyber 205.

The package is organized as follows: Section 0 contains (slightly
amended) {\it netlib} routines, namely {\tt r1mach} and {\tt d1mach},
providing basic machine constants for a variety of computers.
Section 1 contains all the driver routines --- named {\tt test1},
{\tt test2}, etc. --- which are used (and described in the body of this
paper) to test the subroutines of the package.
The complete output of each test is listed immediately after the
driver.
Sections 2--6 constitute the core of the package:
the single- and double-precision subroutines described in the
equally-numbered sections of this paper.
All single-precision routines are provided with comments and
instructions for their use. These, of course, apply to the
double-precision routines as well.

\vspace{.375in}

\noindent
{\small 2.  CLASSICAL WEIGHT FUNCTIONS}

\bigskip

Among the most frequently used orthogonal polynomials are the Jacobi
polynomials, generalized Laguerre polynomials, and Hermite polynomials,
supported respectively on a finite interval, half-infinite interval,
and the whole real line.
The respective weight functions are
$$
w(t) = w^{ ( \alpha , \beta )} (t) =
(1 - t)^{\alpha} (1 + t)^{\beta} ~~\mbox{on}~~
(-1,1), ~~ \alpha > -1, ~ \beta > -1 : ~~
\mbox{Jacobi,}
\eqno(2.1)
$$
$$
w(t) = w^{ ( \alpha )} (t) = t^{\alpha} e^{-t} ~~ \mbox{on} ~~
(0, \infty ) , ~~
\alpha > -1: ~~ \mbox{generalized~Laguerre,}
\eqno(2.2)
$$
$$
w(t) = e^{-t^{2}} ~~ \mbox{on} ~~ (- \infty , \infty ): ~~
\mbox{Hermite.}
\eqno(2.3)
$$
Special cases of the Jacobi polynomials are the Legendre
polynomials ($\alpha = \beta = 0$), the Chebyshev polynomials of the
first ($\alpha = \beta = - \small{\frac{1}{2}}$), second
($\alpha = \beta = \small{\frac{1}{2}}$), third
($\alpha = - \beta = - \small{\frac{1}{2}}$) and fourth
($\alpha = - \beta = \small{\frac{1}{2}}$) kind, and the Gegenbauer polynomials
($\alpha = \beta = \lambda - \small{\frac{1}{2}}$).
The Laguerre polynomials are the special case $\alpha = 0$ of the
generalized Laguerre polynomials.

For each of these polynomials, the corresponding recursion
coefficients $\alpha_k = \alpha_k (w)$,
$\beta_k = \beta_k (w)$ are explicitly known
(see, e.g., [5, pp. 217--221]),
and are generated
in single precision by the routine {\tt recur}.
Its calling sequence is
\smallskip

\begin{center}
{\tt recur(n,ipoly,al,be,a,b,ierr).}
\end{center}
\smallskip

\noindent
\begin{tabular}{lp{3in}}
On entry,  & \\
\\
\ \ \ \ \ {\tt n} & is the number of recursion coefficients
desired; type integer \\
\\
\ \ \ \ \ {\tt ipoly} & is an integer identifying the polynomial
as follows: \\
& 1 = Legendre polynomial on ($-$1,1) \\
& 2 = Legendre polynomial on (0,1) \\
& 3 = Chebyshev polynomial of the first kind \\
& 4 = Chebyshev polynomial of the second kind \\
& 5 = Chebyshev polynomial of the third kind \\
& 6 = Jacobi polynomial with parameters {\tt al,be} \\
& 7 = generalized Laguerre polynomial with pa- \\
& \hspace{.25in} rameter {\tt al} \\
& 8 = Hermite polynomial  \\
\\
\ \ \ \ \ {\tt al,be} & are the input parameters $\alpha , \beta$ for Jacobi and
generalized Laguerre polynomials; type real;  they are
only used if {\tt ipoly} = 6 or 7, and in the latter case
only {\tt al} is used. \\
\end{tabular}

\noindent
\begin{tabular}{lp{3in}}
On return, & \\
\\
\ \ \ \ \ {\tt a,b} & are real arrays of dimension {\tt n} with
{\tt a}($k$), {\tt b}($k$) containing the coefficients
$\alpha_{k-1}$, $\beta_{k-1}$, respectively, $k = 1, 2 ,\ldots,$ {\tt n}
\\
\ \ \ \ \ {\tt ierr} & is an error flag, where \\
& {\tt ierr} = 0 on normal return, \\
& {\tt ierr} = 1 if either {\tt al} or {\tt be} are out of range \\
& \hspace{.5in} when {\tt ipoly} = 6 or {\tt ipoly}=7, \\
& {\tt ierr} = 2 if there is potential overflow in the  \\
&\hspace{.5in} evaluation of $\beta_0$ when {\tt ipoly} = 6 or \\
&\hspace{.5in} {\tt ipoly} = 7;  in this
case $\beta_0$ is set equal \\
& \hspace{.5in} to the largest machine-representable \\
&\hspace{.5in} number, \\
& {\tt ierr} = 3 if {\tt n} is out of range, \\
& {\tt ierr} = 4 if {\tt ipoly} is not one of the \\
&\hspace{.5in} admissible integers.
\end{tabular}
\bigskip

\noindent
No provision has been made for Chebyshev polynomials of the fourth
kind, since their recursion coefficients are obtained from those for
the third-kind Chebyshev polynomials simply by changing the sign
of the $\alpha$'s (and leaving the $\beta$'s unchanged).

The corresponding double-precision routine is {\tt drecur};
it has the same calling sequence, except for real data types now being
double-precision.

In the cases of Jacobi polynomials ({\tt ipoly} = 6) and generalized
Laguerre polynomials ({\tt ipoly} = 7), the recursion coefficient
$\beta_0$ (and only this one) involves
the gamma function $ \Gamma $.
Accordingly, a function routine, {\tt alga}, is provided that computes
the logarithm ln $ \Gamma $
of the gamma function, and a separate routine, {\tt gamma},
computing the gamma function by exponentiating its logarithm.
Their calling sequences are
\bigskip

\hspace{1in}{\tt function alga(x)}

\hspace{1in}{\tt function gamma(x,ierr)},
\medskip

\noindent
where {\tt ierr} is an output variable set equal to
2 or 0 depending on whether the gamma function does, or does not, overflow,
respectively.
The corresponding double-precision routines have the names
{\tt dlga} and {\tt dgamma}.
All these routines require machine-dependent constants
for reasons explained below.

The routine {\tt alga} is based on a rational approximation valid on the
interval $\left[ \small{\frac{1}{2}} , \small{\frac{3}{2}} \right]$.
Outside this interval, the argument $x$ is written as
$$
x = x_e + m,
$$
where
\renewcommand{\arraystretch}{2}
$$
x_e = \left\{ \begin{array}{ll}
x - \lfloor x \rfloor + 1 & \mbox{if} ~~ x - \lfloor x \rfloor \leq \small{\frac{1}{2}} , \\
x - \lfloor x \rfloor & \mbox{otherwise}
\end{array} \right.
$$
\renewcommand{\arraystretch}{1}is in the interval
$\left( \small{\frac{1}{2}} , \small{\frac{3}{2}} \right]$ and $m \geq -1$ is an integer.
If $m = -1$ (i.e., $0 < x \leq \small{\frac{1}{2}}$), then
ln $ \Gamma (x)$ = ln $ \Gamma (x_e ) -$ ln $x$, while for $m > 0$, one
computes ln $ \Gamma (x)$ =\linebreak ln $\Gamma (x_e )$ + ln $p$, where
$p = x_e (x_e + 1) \cdots ( x_e + m - 1 )$.
If $m$ is so large, say $m \geq m_0$, that the product $p$ would overflow, then ln $p$ is computed
(at a price!) as ln $p$ = ln $x_e$ + ln $(x_e +1) + \cdots +$ ln
$(x_e + m-1 )$.
It is here, where a machine-dependent integer is required, namely $m_0$ =
smallest integer $m$ such that $1 \cdot 3 \cdot 5 \cdots (2m+1) / 2^m$ is greater
than or equal to the largest machine-representable number, $R$.
By Stirling's formula, the integer $m_0$ is seen to be the smallest
integer $m$ satisfying $((m+1)/e){\rm ln}((m+1)/e) \geq ({\rm ln} R- \frac{1}{2}
{\rm ln} 8)/e$, hence equal to $ \lfloor e \cdot t(({\rm ln} R- \frac{1}{2} 
{\rm ln} 8)/e) \rfloor$, where $t(y)$ is the inverse function of $y=t \,{\rm ln} t$.
For our purposes, the low-accuracy approximation of $t(y)$,
given in [16, pp. 51--52], and implemented in the routine {\tt t}, is
adequate.

The rational approximation chosen on $\left[ \small{\frac{1}{2}}, \small{\frac{3}{2}} \right]$ is one due to W.J. Cody
and K.E. Hillstrom, namely the one labeled $n = 7$ in Table II of [9].
It is designed to yield
about 16 correct decimal digits (cf. [9, Table I]),
but because of numerical cancellation furnishes only about 13-14
correct decimal digits.

Since rational approximations for ln $ \Gamma $ having sufficient
accuracies for double-precision computation do not seem to be available
in the literature, we used a different approach for the routine
{\tt dlga}, namely the asymptotic approximation (cf. [1, Eq. 6.1.42], where
the constants $B_{2m}$ are Bernoulli numbers)
\renewcommand{\arraystretch}{2}
$$
\begin{array}{lll}
\mbox{ln} ~\Gamma (y) & = & \left( y - \small{\frac{1}{2}} \right)
\mbox{ln} ~ y - y + \small{\frac{1}{2}} ~\mbox{ln} ~(2 \pi ) \\
& + & \sum_{m=1}^n ~\frac{B_{2m}}{2m(2m-1)} ~y^{-(2m-1)} + R_n (y)
\end{array}
\eqno(2.4)
$$
\renewcommand{\arraystretch}{1}for values of
$y > 0$ large enough to have
$$
| R_n (y) | \leq ~\frac{1}{2}~ 10^{-d} ,
\eqno(2.5)
$$
where $d$ is the number of decimal digits carried in double-precision
arithmetic, another machine-dependent real number.
If (2.5) holds for $y \geq y_0$,
and if $x \geq y_0$, we compute ln $\Gamma (x)$ from the asymptotic
expression (2.4) (where $y = x$
and the remainder term is neglected).
Otherwise, we let
$k_0$ be  the smallest positive integer $k$
such that $x + k \geq y_0$, and use
$$
\mbox{ln} ~\Gamma (x) = \mbox{ln} ~\Gamma (x + k_0 ) -
\mbox{ln} ~(x (x+1 ) \cdots (x + k_0 - 1) ) ,
\eqno(2.6)
$$
where the first
term on the right is computed
from (2.4) (with $y = x + k_0$).
Since for $y > 0$,
$$
| R_n (y) | \leq ~
\frac{| B_{2n+2}|}{(2n+2)(2n+1)} ~y^{-(2n+1)}
$$
(cf. [1, Eq. 6.1.42]), the inequality (2.5) is satisfied if
$$
y \geq \mbox{exp} \left\{
\frac{1}{2n+1} ~\left[ d~ \mbox{ln}~ 10 + \mbox{ln} ~
\frac{2 | B_{2n+2} |}{(2n+2)(2n+1)} ~\right] \right\} .
\eqno(2.7)
$$
In our routine {\tt dlga}, we have selected $n = 9$.
For double-precision accuracy on the Cyber 205, we have
$d \approx 28.3$, for which (2.7) then gives
$y \geq \mbox{exp} \{ .121188\ldots ~ d + .053905\ldots \} \approx 32.6$.

For single-precision calculation we selected the method of rational
approximation, rather than the asymptotic formula (2.4) and (2.6),
since we found that the former is generally more accurate, by a factor, on
the average, of about 20 and as large as 300.
Neither method yields full machine accuracy.
The former, as already mentioned, loses accuracy in the evaluation of the
approximation.
The latter suffers loss of accuracy because of cancellation
occurring in (2.6), which typically amounts to a loss of 2--5 significant decimal digits
in the gamma function itself.

Since these errors affect only the coefficient $\beta_0$ (and only
if {\tt ipoly} = 6 or 7), they are of no consequence unless the output
of the routine {\tt recur} serves as input to another routine,
such as {\tt gauss} (cf. \S6), which makes essential use of $\beta_0$.
In this case, for maximum single-precision accuracy, it is recommended that $\beta_0$
be first obtained in double precision by means of {\tt drecur} with
{\tt n} = 1, and then converted to single precision.

\vspace{.375in}

\noindent
{\small 3.  MOMENT-RELATED METHODS}
\bigskip

It is a well-known fact that the first $n$ recursion coefficients
$\alpha_k ( d \lambda )$,
$\beta_k ( d \lambda )$, $k = 0, 1, \ldots , n-1$ (cf. (1.3)), are uniquely
determined by the first $2n$ moments $\mu_k$
of the measure $d \lambda$,
$$
\mu_k = \mu_k ( d \lambda ) = \int_{\mbox{\myfont R}}
t^k d \lambda (t) , ~~~~
k = 0, 1, 2 ,\ldots, 2n-1 .
\eqno(3.1)
$$
Formulae are known, for example, which express the
$\alpha$'s and $\beta$'s in terms of Hankel determinants in these
moments.
The trouble is that these formulae become increasingly sensitive to small
errors as $n$ becomes large.
There is an inherent reason for this:
the underlying (nonlinear) map
$K_n$: $\mbox{\myfont R}^{2n} \rightarrow \mbox{\myfont R}^{2n}$ has
a condition number, cond $K_n$, which grows exponentially with $n$
(cf. [19, \S3.2]).
Any method that attempts to compute the desired coefficients from the
moments in (3.1), therefore, is doomed to fail, unless $n$ is quite
small, or extended precision is being employed.
That goes in particular for an otherwise elegant method due to Chebyshev
(who developed it for the case of discrete measures $d \lambda$) that
generates the $\alpha$'s and $\beta$'s directly from the moments (3.1),
bypassing determinants altogether (cf. [4], [19,\S2.3]).

Variants of Chebyshev's algorithm with more satisfactory stability
properties
have been developed
by Sack and Donovan [46] and Wheeler [50]
(independently of Chebyshev's work).
The idea is to forego the moments (3.1) as input data, and instead
depart from so-called {\it modified moments}.
These are defined by replacing the power $t^k$ in (3.1) by an
appropriate polynomial $p_k (t)$ of degree $k$,
$$
\nu_k = \nu_k ( d \lambda ) = \int_{\mbox{\myfont R}}
p_k (t) d \lambda (t) , ~~~~
k = 0, 1, 2 ,\ldots, 2n-1 .
\eqno(3.2)
$$
For example, $p_k$ could be one of the classical orthogonal polynomials.
More generally, we shall assume that $\{ p_k \}$ are monic
polynomials satisfying a three-term recurrence relation similar to the
one in (1.3),
\renewcommand{\arraystretch}{2}
$$
\begin{array}{c}
{\displaystyle
p_{k+1} (t) = (t - a_k ) p_k (t) - b_k p_{k-1} (t) , ~~~~
k = 0, 1, 2,\ldots, } \\
{\displaystyle
p_0 (t) = 1 , ~~~~ p_{-1} (t) = 0 , }
\end{array}
\eqno(3.3)
$$
\renewcommand{\arraystretch}{1}with coefficients
$a_k \in \mbox{\myfont R}$, $b_k \geq 0$ that are known.
(In the special case
$a_k = 0$, $b_k = 0$, we are led back to powers and
ordinary moments.)
There now exists an algorithm, called {\it modified Chebyshev algorithm}
in [19, \S2.4], which takes the $2n$ modified moments in (3.2) and the
$2n-1$ coefficients $\{ a_k \}_{k=0}^{2n-2}$,
$\{ b_k \}_{k=0}^{2n-2}$ in (3.3), and from them generates the $n$
desired coefficients $\alpha_k ( d \lambda )$,
$\beta_k ( d \lambda )$, $k = 0, 1,\ldots, n-1$.
It generalizes Chebyshev's algorithm, which
can be recovered (if need be) by putting
$a_k = b_k = 0$.
The modified Chebyshev algorithm is embodied in the subroutine
{\tt cheb}, which has the calling sequence
\medskip

\noindent
\hspace{.75in} {\tt cheb(n,a,b,fnu,alpha,beta,s,ierr,s0,s1,s2)}

\noindent
\hspace{.75in} {\tt dimension a(*),b(*),fnu(*),alpha(n),beta(n),s(n),}

\noindent
\hspace{.75in} {\tt \ \ s0(*),s1(*),s2(*)}

\noindent
\begin{tabular}{lp{3in}}
\\
On entry,  & \\
\\
\ \ \ \ \ {\tt n} & is the number of recursion coefficients
desired; type integer \\
\\
\ \ \ \ \ {\tt a,b} & are arrays of dimension 2$\times${\tt n}--1 holding
the coefficients {\tt a}$(k) = a_{k-1}$, {\tt b}$(k) = b_{k-1}$,
$k = 1, 2,\ldots, 2n-1$ \\ 
\\
\ \ \ \ \ {\tt fnu} & is an array of dimension 2$\times${\tt n} holding the modified
moments {\tt fnu}$(k) = \nu_{k-1}$, $k = 1, 2,\ldots,2\times{\tt n}$ \\
\\
On return, & \\
\\
\ \ \ \ \ {\tt alpha,beta} & are arrays of dimension {\tt n} containing
the desired recursion coefficients {\tt alpha} $(k) = \alpha_{k-1}$,
{\tt beta} $(k) = \beta_{k-1}$, $k = 1, 2,\ldots,$ {\tt n} \\
\\
\ \ \ \ \ {\tt s} & is an array of dimension {\tt n} containing
the numbers {\tt s}$(k) = \int_{\mbox{\myfont R}} \pi_{k-1}^2 d \lambda$,
$k = 1, 2, \ldots,$ {\tt n} \\
\\
\ \ \ \ \ {\tt ierr} & is an error flag, equal to 0 on normal return,
equal to 1 if $|\nu_0|$ is less than the machine zero, equal to 2 if
{\tt n} is out of range, equal to
$-(k-1)$ if {\tt s}$(k)$, $k = 1, 2,\ldots,$ {\tt n}, is about
to underflow, and equal to $+(k-1)$, if it is about to overflow.
\end{tabular}
\bigskip

\noindent
The arrays {\tt s0}, {\tt s1}, {\tt s2} of dimension 2$\times${\tt n}
are needed for working space.
\bigskip

There is again a map
$K_n$: $\mbox{\myfont R}^{2n} \rightarrow \mbox{\myfont R}^{2n}$
underlying the modified Chebyshev algorithm, namely the map taking the
2$n$ modified moments into the $n$ pairs of recursion coefficients.
The condition of the map $K_n$ (actually of a somewhat different, but
closely related map) has been studied in [19, \S3.3] and [23] in the
important case where the polynomials $p_k$ defining the modified moments
are themselves orthogonal polynomials,
$p_k ( \cdot ) = p_k (~ \cdot ~; d \mu )$, with respect to a measure
$d \mu$ (for example, one of the classical ones) for which the recursion
coefficients $a_k$, $b_k$ are known.
The upshot of the analysis then is that the condition of $K_n$ is
characterized by a certain positive polynomial
$g_n (~\cdot~ ; d \lambda )$ of degree $4n-2$, depending only on the target
measure $d \lambda$, in the sense that
$$
\mbox{cond}~K_n = \int_{\mbox{\myfont R} } g_n
(t ; d \lambda ) d \mu (t) .
\eqno(3.4)
$$
Thus, the numerical stability of the modified Chebyshev algorithm
is determined by the magnitude of $g_n$ on the support of $d \mu$.

The occurrence of underflow [overflow] in the computation of the
$\alpha$'s and $\beta$'s, especially on computers with limited
exponent range, can often be avoided by multiplying all modified
moments by a sufficiently large [small] scaling factor before entering
the routine. On exit, the coefficient $\beta_{0}$ (and only this one!)
then has to be divided by the same scaling factor. (There may occur
harmless underflow of auxiliary quantities in the routine {\tt cheb},
which is difficult to avoid since some of these quantitites actually
are expected to be zero.)
\bigskip

{\it Example} 3.1.
$d \lambda_\omega (t) = [(1 - \omega^2 t^2 ) (1-t^2 )]^{- 1/2} dt$ on
($-$1,1), $0 \leq \omega < 1$.
\smallskip

This example is of some historical interest, in that it has
already been considered by Christoffel [8, Example 6];
see also [44].
Computationally, the example is of interest as there are empirical reasons to
believe that for the choice
$d \mu (t) = (1-t^2 )^{- 1/2} dt$ on ($-$1,1) --- which appears rather natural
--- the modified Chebyshev algorithm is exceptionally stable, uniformly in $n$,
in the sense that in (3.4) one has $g_n \leq 1$ on
supp $d \mu$ for all $n$ (cf. [22, Example 5.2]).
With the above choice of $d \mu$, the polynomials $p_k$ are clearly the
Chebyshev polynomials of the first kind,
$p_0 = T_0$,
$p_k = 2^{-(k-1)} T_k$,
$k \geq 1$, and the modified moments are given by
$$
\nu_0 = \int_{-1}^1 d \lambda_{\omega} (t) , ~~
\nu_k = ~\frac{1}{2^{k-1} } ~
\int_{-1}^1 T_k (t) d \lambda_{\omega} (t) , ~~
k = 1, 2, 3,\ldots ~.
\eqno(3.5)
$$
They are expressible in terms of the Fourier coefficients
$C_r ( \omega^2 )$ in
$$
(1 - \omega^2 ~\sin^2 \theta)^{- 1/2} =
C_0 ( \omega^2 ) + 2
\sum_{r=1}^\infty C_r ( \omega^2 ) \cos 2 r \theta
\eqno(3.6)
$$
by means of (cf. [19, Example 3.3])
\renewcommand{\arraystretch}{2}
$$
\begin{array}{l} ~~
\nu_0 = \pi C_0 ( \omega^2 ), \\
\left.
\begin{array}{l}
\nu_{2m} = ~ \frac{(-1)^m \pi}{2^{2m-1}} ~C_m ( \omega^2 ) \\
\nu_{2m-1} = 0 ~~~~~~~~~~
\end{array}
\right\} ~~~
m = 1, 2, 3, \ldots ~.
\end{array}
\eqno(3.7)
$$
\renewcommand{\arraystretch}{1}The Fourier
coefficients
$\{ C_r ( \omega^2 ) \}$ in turn
can be accurately computed as minimal solution of a certain
three-term recurrence relation (see [19, pp. 310--311]).

The ordinary moments
$$
\mu_0 = \nu_0 , ~~
\mu_k = \int_{-1}^1 t^k d \lambda_{\omega} (t) , ~~
k = 1, 2, 3,\ldots,
\eqno(3.8)
$$
likewise
can be expressed in terms of the Fourier coefficients
$C_r ( \omega^2 )$ by writing $t^{2m}$ as a linear combination of
Chebyshev polynomials $T_0$, $T_2 ,\ldots, T_{2m}$ (cf. Eq. (22) on
p. 454 of [43]).
The result is
\renewcommand{\arraystretch}{2}
$$
\begin{array}{l}
{\displaystyle
\mu_{2m} = ~\frac{(-1)^m \pi}{2^{2m-1}} ~
\sum_{r=0}^m (-1)^r \gamma_r^{(m)} C_{m-r} ( \omega^2) ,} \\
{\displaystyle
\mu_{2m-1} = 0 ,} ~~~
\end{array}
m = 1, 2, 3,\ldots,
\eqno(3.9)
$$
\renewcommand{\arraystretch}{1}where
\renewcommand{\arraystretch}{2}
$$
\begin{array}{l}
{\displaystyle
\gamma_0^{(m)} = 1 ,} \\
{\displaystyle
\gamma_r^{(m)} = ~\frac{2m+1-r}{r}~ \gamma_{r-1}^{(m)} , ~~
r = 1, 2,\ldots, m-1,} \\
{\displaystyle
\gamma_m^{(m)} = ~\frac{m+1}{2m} ~\gamma_{m-1}^{(m)}.}
\end{array}
\eqno(3.10)
$$
\renewcommand{\arraystretch}{1}

The driver routine {\tt test1} (in \S1 of the package) generates for
$\omega^2$ = .1(.2).9, .99, .999 the first $n$ recurrence coefficients
$\beta_k (d \lambda_{\omega} )$ (all $\alpha_k = 0$), both in single and double precision,
using modified moments if {\tt modmom}={\tt .true.}, and ordinary
moments otherwise.
In the former case, $n$ = 80, in the latter, $n$ = 20.
It prints the double-precision values of $\beta_k$ together with the
relative errors of the single-precision values
(computed as the difference of the double-precision and single-precision
values divided by the double-precision value).
In {\tt test1}, as well as in all subsequent drivers, not all error
flags are interrogated for possible malfunction.
The user is urged, however, to do so as a matter of principle.

\bigskip
\begin{center}
\begin{tabular}{llll}
\multicolumn{1}{c}{$\omega^2$} &
\multicolumn{1}{c}{$k$} &
\multicolumn{1}{c}{$\beta_k^{~\rm double}$} &
\multicolumn{1}{c}{$\mbox{err}~ \beta_k^{~\rm single}$} \\ \hline
 & & &  \\*[-10pt]
.1 & 0 & 3.224882697440438796459832725 & 1.433(--14) \\
& 1 & \ \ .5065840806382684475158495727 & 1.187(--14) \\
& 5 & \ \ .2499999953890031901881028267 & 1.109(--14) \\
& 11 & \ \ .2499999999999999996365048540 & 1.454(--18) \\
& 18 &\ \ .2500000000000000000000000000 & 0.000 \\
.5 & 0 & 3.708149354602743836867700694 & 9.005(--15) \\
& 1 & \ \ .5430534189555363746250333773 & 2.431(--14) \\
& 8 & \ \ .2499999846431723296083779480 & 4.109(--15) \\
& 20 & \ \ .2499999999999999978894635584 & 8.442(--18) \\
& 35 & \ \ .2500000000000000000000000000 & 0.000 \\
.9 & 0 & 5.156184226696346376405141543 & 6.950(--15) \\
& 1 & \ \ .6349731661452458711622492613 & 7.920(--15) \\
& 19 & \ \ .2499999956925950094629502830 & 1.820(--14) \\
& 43 & \ \ .2499999999999998282104100896 & 6.872(--16) \\
& 79 & \ \ .2499999999999999999999999962 & 1.525(--26) \\
.999 & 0 & 9.682265121100594060678208257 & 1.194(--13) \\
& 1 & \ \ .7937821421385176965531719571 & 6.311(--14) \\
& 19 & \ \ .2499063894398209200047452537 & 1.026(--14) \\
& 43 & \ \ .2499955822633680825859750068 & 8.282(--15) \\
& 79 & \ \ .2499998417688157876153069211 & 1.548(--15)
\end{tabular}
\end{center}
\medskip
\begin{center}
TABLE 3.1.  Selected output from {\tt test1} in the case \\
of modified moments
\end{center}
\bigskip

The routine
\smallskip

\begin{center}
{\tt fmm(n,eps,modmom,om2,fnu,ierr,f,f0,rr)}
\end{center}
\smallskip

\noindent
used by the driver computes the first 2$\times${\tt n} modified [ordinary] moments
for
$\omega^2$ = {\tt om}$2$, to a relative accuracy {\tt eps}
if {\tt modmom}={\tt .true.} [{\tt .false.}].
The results are stored in the array {\tt fnu}.
The arrays {\tt f}, {\tt f0}, {\tt rr} are internal working arrays of
dimension {\tt n}, and {\tt ierr} is an error flag.
On normal return, {\tt ierr} = 0; otherwise,
{\tt ierr} = 1, indicating lack of convergence (within a prescribed
number of iterations) of the backward recurrence algorithm for computing the
minimal solution $\{ C_r ( \omega^2 ) \}$.
The latter is likely to occur if $\omega^2$ is too close to 1.
The routine {\tt fmm} as well as its double-precision version
{\tt dmm} are listed
immediately after the routine {\tt test1}.

In Table 3.1, we show
selected results from the output of {\tt test1}, when
{\tt modmom}={\tt .true.} (Complete results are given in the
package immediately after {\tt test1}.)
The values for $k = 0$ are expressible in terms of the complete elliptic
integral,
$\beta_0 = 2 K ( \omega^2 )$, and were checked, where
possible, against the 16S-values in Table 17.1 of [1].
In all cases,
there was agreement to all 16 digits.
The largest relative error observed was $2.43 \times 10^{-13}$ for
$\omega^2 = .999$ and $k = 2$.
When $\omega^2 \leq .99$,
the error was always less than
$2.64 \times 10^{-14}$, which confirms the extreme stability
of the modified Chebyshev algorithm in this example.
It can be seen (as was to be expected) that for $\omega^2$
not too close to 1, the coefficients converge rapidly to
$\small{\frac{1}{4}}$.

In contrast, Table 3.2 shows selected results (for complete results,
see the package) in the case of
ordinary moments ({\tt modmom}={\tt .false.}) and demonstrates
the severe instability of the Chebyshev algorithm.
We remark that the moments themselves are all accurate to essentially
machine precision, as has been verified by additional computations.

\bigskip
\begin{center}
\begin{tabular}{rrl|lrl}
\multicolumn{1}{c}{$\omega^2$} &
\multicolumn{1}{c}{$k$} &
\multicolumn{1}{c|}{err $\beta_k$} &
\multicolumn{1}{c}{$\omega^2$} &
\multicolumn{1}{c}{$k$} &
\multicolumn{1}{c}{err $\beta_k$} \\ \hline
 & & & & &  \\*[-10pt]
.1 & 1 & 1.187(--14) & .9 & 1 & 3.270(--15) \\
& 7 & 2.603(--10) & & 7 & 4.819(--10) \\
& 13 & 9.663(--6) & & 13 & 1.841(--5) \\
& 19 & 4.251(--1) & & 19 & 6.272(--1) \\
.5 & 1 & 2.431(--14) & .999 & 1 & 6.311(--14) \\
& 7 & 5.571(--10) & & 7 & 1.745(--9) \\
& 13 & 9.307(--6) & & 13 & 8.589(--5) \\
& 19 & 5.798(--1) & & 19 & 4.808(0)
\end{tabular}
\end{center}
\medskip
\begin{center}
TABLE 3.2.  Selected output from {\tt test1} in the case of \\
ordinary moments
\end{center}
\bigskip

The next example deals with another weight function for which the modified
Chebyshev algorithm performs rather well.

\bigskip

{\it Example} 3.2.
$d \lambda_\sigma (t) = t^{\sigma}~\mbox{ln} ~(1/t) dt$ on (0,1], $\sigma > -1$.
\smallskip

What is nice about this example is that both modified and ordinary moments
of $d \lambda_{\sigma}$ are known in closed form.
The latter are obviously given by
$$
\mu_k (d \lambda_{\sigma} ) = ~\frac{1}{\sigma + 1 + k} ~, ~~~~
k = 0, 1, 2 ,\ldots,
\eqno(3.11)
$$
whereas the former, relative to shifted monic Legendre polynomials
({\tt ipoly}=2 in {\tt recur}), are (cf. [17])
\renewcommand{\arraystretch}{2}
$$
\frac{(2k)!}{k!^2} \nu_k (d \lambda_{\sigma} ) = \left\{
\begin{array}{l}
{\displaystyle
(-1)^{k - \sigma} ~\frac{\sigma !^2 (k - \sigma -1)!}{(k+ \sigma +1)!} ~,~~~~
0 \leq \sigma < k, ~~\sigma~\mbox{an~integer},} \\
{\displaystyle
\frac{1}{\sigma +1} ~\left\{ \frac{1}{\sigma +1} ~+
\sum_{r=1}^k \left( \frac{1}{\sigma +1+r} ~-~
\frac{1}{\sigma + 1-r} \right) \right\}
\prod_{r=1}^k ~\frac{\sigma +1-r}{\sigma +1+r} ~,} \\
{\displaystyle
\hspace{2in} \mbox{otherwise}.}
\end{array}
\right.
\eqno(3.12)
$$
\renewcommand{\arraystretch}{1}The routines
{\tt fmm} and {\tt dmm} appended to {\tt test2} in \S1 of the
package, similarly as the corresponding
routines in Example 3.1, generate the first
2$\times${\tt n} modified moments
$\nu_0 , \nu_1 ,\ldots, \nu_{2n-1}$, if {\tt modmom}={\tt .true.},
and the first 2$\times${\tt n} ordinary moments, otherwise.
The calling sequence of {\tt fmm} is
\smallskip

\begin{center}
{\tt fmm(n,modmom,intexp,sigma,fnu).}
\end{center}
\smallskip

\noindent
The logical variable {\tt intexp} is to be set {\tt .true.}, if
$\sigma$ is an integer, and {\tt .false.} otherwise.
In either case, the input variable {\tt sigma} is assumed of type
{\tt real}.

The routine {\tt test2} generates the first {\tt n} recursion
coefficients $\alpha_k ( d \lambda_{\sigma} )$,
$\beta_k (d \lambda_{\sigma} )$ in single and double precision
for $\sigma = - \small{\frac{1}{2}}, 0, \small{\frac{1}{2}}$,
where {\tt n} = 100 for the modified Chebyshev algorithm
({\tt modmom=.true.}), and {\tt n} = 12
for the classical Chebyshev algorithm ({\tt modmom=.false.}).
Selected double-precision results to 25 significant digits,
when modified moments are used,
are shown in Table 3.3. (The complete results are given in the
package after {\tt test2}.)

\bigskip
\begin{center}
\begin{tabular}{rrrl}
\multicolumn{1}{c}{$\sigma$} &
\multicolumn{1}{c}{$k$} &
\multicolumn{1}{c}{$\alpha_k$} &
\multicolumn{1}{c}{$\beta_k$} \\ \hline
\\
--.5 & 0 & .1111111111111111111111111 & 4.000000000000000000000000 \\
& 12 & .4994971916094638566242202 & \ \ .06231277082877488477563886 \\
& 24 & .4998662912324218943801592 & \ \ .06245372557342242600457226 \\
& 48 & .4999652635485445800661969 & \ \ .06248855717748684742433618 \\
& 99 & .4999916184024356271670789 & \ \ .06249733823051821636937156 \\
0 & 0  & .2500000000000000000000000 & 1.000000000000000000000000  \\
& 12 & .4992831802157361310272625 & \ \ .06238356835953571123560330 \\
& 24 & .4998062839486146398501532 & \ \ .06247100084469111001639128 \\
& 48 & .4999494083797023879356424 & \ \ .06249281268110967462373889 \\
& 99 & .4999877992015903283047919 & \ \ .06249832670616925926204896 \\
.5 & 0 & .3600000000000000000000000 & \ \ .4444444444444444444444444 \\
& 12 & .4993755732917555644203267 & \ \ .06237082738280752611960887 \\
& 24 & .4998324497706394488722725 & \ \ .06246581011945496883543089 \\
& 48 & .4999567275223771727791521 & \ \ .06249115332711027176695932 \\
& 99 & .4999896931841789781887674 & \ \ .06249787251281682973825635
\end{tabular}
\end{center}
\medskip
\begin{center}
TABLE 3.3.  Selected output from {\tt test2} in the case of \\
modified moments
\end{center}
\bigskip

\noindent
The largest relative errors observed, over all $k = 0, 1,\ldots, 99$,
were respectively $6.211 \times 10^{-11}$,
$2.237 \times 10^{-12}$,
$1.370 \times 10^{-12}$ for the $\alpha$'s, and
$1.235 \times 10^{-10}$,
$4.446 \times 10^{-12}$,
$2.724 \times 10^{-12}$ for
the $\beta$'s, attained consistently at $k$ = 99.
The accuracy achieved is slightly less than in Example 3.1, for
reasons explained in [22, Example 5.3].

The complete results for $\sigma =~-~\small{\frac{1}{2}}$ are also
available in [27, Appendix, Table 1].
(They differ occasionally by 1 unit in the last decimal place from those
produced here, probably because of a slightly different computation of the
modified moments.)
The results for $\sigma = 0$ can be checked up to $k = 15$
against the 30S-values given in [47, p. 92],
and for $16 \leq k \leq 19$ against 12S-values in [10,
Table 3].
There is complete agreement to all 25 digits in the former
case, and agreement to 12 digits in the latter, although there are
occasional end figure discrepancies of 1 unit.
These are believed to be due to rounding errors committed in [10],
since similar discrepancies occur also in the range $k \leq 15$.
We do not know of any tables for $\sigma = ~\small{\frac{1}{2}}$,
but a test will be given in \S5, Example 5.1.

The use of ordinary moments ({\tt modmom}= {\tt .false.}) produces predictably worse
results, the relative errors of which are shown in Table 3.4.

\bigskip
\begin{center}
\begin{tabular}{rlllllllll}
\multicolumn{1}{c}{$k$} &
\multicolumn{1}{c}{$\sigma$} &
\multicolumn{1}{c}{err $\alpha_k$} &
\multicolumn{1}{c}{err $\beta_k$} &
\multicolumn{1}{c}{$\sigma$} &
\multicolumn{1}{c}{err $\alpha_k$} &
\multicolumn{1}{c}{err $\beta_k$} &
\multicolumn{1}{c}{$\sigma$} &
\multicolumn{1}{c}{err $\alpha_k$} &
\multicolumn{1}{c}{err $\beta_k$} \\ \hline
 & & & & & & & & &  \\*[-10pt]
2 & --.5 & 1.8(--13) & 7.7(--14) & 0 & 4.2(--13) & 7.6(--13) & .5 & 1.6(--12) & 2.6(--13) \\
5 & & 2.2(--9) & 1.2(--9) & & 4.2(--9) & 1.2(--10) & & 1.3(--8) & 6.6(--9) \\
8 & & 1.1(--5) & 5.5(--6) & & 4.3(--6) & 3.8(--6) & & 6.0(--5) & 5.2(--6) \\
11 & & 2.5(--1) & 1.7(--1) & & 1.3(0) & 3.2(--1) & & 2.2(0) & 4.7(--1)
\end{tabular}
\end{center}
\medskip
\begin{center}
TABLE 3.4.  Selected output from {\tt test2} in the case of \\
ordinary moments
\end{center}

\vspace{.375in}

\noindent
{\small 4.  STIELTJES, ORTHOGONAL REDUCTION, AND DISCRETIZATION PROCEDURES}
\bigskip

4.1. {\it The Stieltjes procedure}.
It is well known that the coefficients $\alpha_k ( d \lambda )$,
$\beta_k ( d \lambda )$ in the basic recurrence relation (1.3) can be
expressed in terms of the orthogonal polynomials (1.2) and the inner
product (1.1) as follows:
\renewcommand{\arraystretch}{2}
$$
\begin{array}{c}
{\displaystyle
\alpha_k (d \lambda ) = ~
\frac{(t \pi_k , \pi_k )}{( \pi_k , \pi_k )} , ~~
k \geq 0 ;} \\
{\displaystyle
\beta_0 (d \lambda ) = ( \pi_0 , \pi_0 ) , ~~
\beta_k (d \lambda ) = ~
\frac{(\pi_k , \pi_k )}{( \pi_{k-1} , \pi_{k-1} )} , ~~
k \geq 1.}
\end{array}
\eqno(4.1)
$$
\renewcommand{\arraystretch}{1}Provided
the inner product can be readily calculated, (4.1) suggests
the following ``bootstrap'' procedure:
Compute $\alpha_0$ and $\beta_0$ by the
first relations in (4.1) for $k = 0$.
Then use the recurrence relation (1.3) for $k = 0$ to obtain $\pi_1$.
With $\pi_0$ and $\pi_1$ known, apply (4.1) for $k = 1$ to get
$\alpha_1$, $\beta_1$, then again (1.3) to obtain $\pi_2$, and so on.
In this way, alternating between (4.1) and (1.3), we can bootstrap
ourselves up to as many of the coefficients $\alpha_k$, $\beta_k$ as are desired.
We attributed this procedure to Stieltjes, and called it
{\it Stieltjes's procedure} in [19].

In the case of discrete orthogonal polynomials, i.e., for inner
products of the form
$$
(u,v) = \sum_{k=1}^N w_k
u(x_k ) v(x_k ) , ~~~~ w_k > 0 ,
\eqno(4.2)
$$
Stieltjes's procedure is easily implemented; the resulting routine is
called {\tt sti}, and has the calling sequence
\medskip

\begin{center}
{\tt
sti(n,ncap,x,w,alpha,beta,ierr,p0,p1,p2)}.
\end{center}
\medskip

\noindent
\begin{tabular}{lp{3in}}
On entry,  & \\
\\
\ \ \ \ \ {\tt n} & is the number of recursion coefficients
desired; type integer \\
\\
\ \ \ \ \ {\tt ncap} & is the number of terms, $N$, in the discrete
inner product; type integer \\
\\
\ \ \ \ \ {\tt x,w} & are arrays of dimension {\tt ncap} holding
the abscissae {\tt x} $(k) = x_k$ and weights {\tt w} $(k) = w_k$,
$k = 1, 2,\ldots,$ {\tt ncap},
of the discrete inner product. \\
\\
On return, & \\
\\
\ \ \ \ \ {\tt alpha,beta} & are arrays of dimension {\tt n} containing
the desired recursion coefficients {\tt alpha} $(k) = \alpha_{k-1}$,
{\tt beta} $(k) = \beta_{k-1}$, $k = 1, 2,\ldots,$ {\tt n} \\
\\
\ \ \ \ \ {\tt ierr} & is an error flag having the value 0 on normal
return, and the value 1 if {\tt n} is not in the proper range
$1 \leq n \leq N$; if during the computation of a recursion coefficient
with index $k$ there is impending underflow or overflow, {\tt ierr} will
have the value $-k$ in case of underflow, and the value $+k$ in case
of overflow. (No error flag is set in case of harmless underflow.)
\end{tabular}
\bigskip

\noindent
The arrays {\tt p0}, {\tt p1}, {\tt p2} are working arrays of dimension
{\tt ncap}.
The double-precision routine has the name {\tt dsti.}

Occurrence of underflow [overflow] can be forestalled by multiplying
all weights $w_{k}$ by a sufficiently large [small] scaling factor prior
to entering the routine. Upon return, the coefficient $\beta _{0}$
will then have to be readjusted by dividing it by the same scaling
factor.
\bigskip

4.2. {\it Orthogonal reduction method}.
Another approach to producing the recursion coefficients
$\alpha_k$, $\beta_k$ from the quantities $x_k$, $w_k$ defining
the inner product (4.2) is based on the observation
(cf. [2], [29, \S7]) that the symmetric tridiagonal matrix of
order $N + 1$,

$$
J (d \lambda_N ) = \left[ \begin{array}{ccccc}
1 & \sqrt{\beta_0} & & & 0 \\
\sqrt{\beta_0} & \alpha_0 & \sqrt{\beta_1} & & \\
& \sqrt{\beta_1} & \alpha_1 & \ddots & \\
& & \ddots & \ddots & \sqrt{\beta_{N-1}} \\
0 & & & \sqrt{\beta_{N-1}} & \alpha_{N-1}
\end{array}
\right]
\eqno(4.3)
$$
(the ``extended Jacobi matrix'' for the discrete measure
$d \lambda_N$ implied in (4.2)), is orthogonally similar to
the matrix
$$
\left[ \begin{array}{cc}
1 & \sqrt{w}^T \\
& \\
\sqrt{w} & D_x
\end{array}
\right] , ~~
\sqrt{w} = \left[ \begin{array}{c}
\sqrt{w_1} \\ \vdots \\ \sqrt{w_N}
\end{array}
\right] , ~~
D_x = \left[ \begin{array}{c}
x_1 ~~~~ 0 \\ \ddots \\ ~~ 0 ~~~~ x_N
\end{array}
\right] .
\eqno(4.4)
$$
Hence, the desired matrix $J(d \lambda_N )$ can be
obtained by applying Lanczos's algorithm to the matrix (4.4).
This is implemented in the routine
\medskip

\begin{center}
{\tt lancz(n,ncap,x,w,alpha,beta,ierr,p0,p1)},
\end{center}
\medskip

\noindent
which uses a judiciously constructed sequence of Givens transformations
to accomplish the orthogonal similarity transformation (cf.
[45,3,41,2]; the routine {\tt lancz} is adapted from the
routine RKPW in [41, p. 328]).
The input and output parameters of the routine {\tt lancz} have the same
meaning as in the routine {\tt sti}, except that {\tt ierr} can only
have the value 0 or 1, while {\tt p0}, {\tt p1} are
again working arrays of dimension {\tt ncap}.
The double-precision version of the routine is named {\tt dlancz}.

The routine {\tt lancz} is generally superior to the routine 
{\tt sti}: The procedure used in {\tt sti}
may develop numerical instability from some point on, and therefore
give inaccurate results for larger values of {\tt n}. It furthermore
is subject to underflow and overflow conditions. None of these
shortcomings is shared by the routine {\tt lancz}. On the other hand,
there are cases where {\tt sti} does better than {\tt lancz} (cf.
Example 4.5).

We illustrate the phenomenon of instability (which is explained in [32])
in the case of the ``discrete Chebyshev'' polynomials.

\bigskip

{\it Example} 4.1.
The inner product (4.2) with
$x_k = -1 + 2 ~\frac{k-1}{N-1}$, $~w_k = ~\frac{2}{N}$, $~k = 1, 2,\ldots, N$.
\smallskip

This generates discrete analogues of the Legendre polynomials, which they indeed
approach as $N \rightarrow \infty$.
The recursion coefficients are explicitly known:

\renewcommand{\arraystretch}{2}
$$
\begin{array}{c}
{\displaystyle
\alpha_k = 0, ~~~~ k = 0, 1,\ldots, N-1;} \\
{\displaystyle
\beta_0 = 2, ~~ \beta_k = \left(
1 + \frac{1}{N-1} \right)^2 \left( 1 - \left( \frac{k}{N} \right)^2 \right) ~~
\left( 4 - \frac{1}{k^2} \right)^{-1} } \\
{\displaystyle
k = 1, 2,\ldots, N-1 . }
\end{array}
\eqno(4.5)
$$
\renewcommand{\arraystretch}{1}To find out how well the
routines {\tt sti} and {\tt lancz} generate them 
(in single precision), when $N$ = 40, 80, 160, 320,
we wrote the driver {\tt test3}, which computes the
respective absolute errors for the $\alpha$'s
and relative errors for the $\beta$'s.

Selected results for Stieltjes's algorithm are shown in Table 4.1.
The gradual deterioration, after some point (depending on $N$), is
clearly visible.
Lanczos's method, in contrast, preserves essentially full accuracy;
the largest error in the $\alpha$'s is
1.42(--13), 2.27(--13), 4.83(--13),
8.74(--13) for $N$ = 40, 80, 160, 320,
respectively, and 3.38(--13), 6.63(--13), 2.17(--12),
5.76(--12) for the $\beta$'s.
\bigskip

\noindent
\begin{tabular}{crllclll}
\multicolumn{1}{c} {$N$} &
\multicolumn{1}{c} {$n$} &
\multicolumn{1}{c} {err $\alpha$} &
\multicolumn{1}{c} {err $\beta$} &
\multicolumn{1}{c} {$N$} &
\multicolumn{1}{c} {$n$} &
\multicolumn{1}{c} {err $\alpha$} &
\multicolumn{1}{c} {err $\beta$} \\ \hline
 & & & & & & &   \\*[-10pt]
40 & $\leq$ 35 & $\leq$ 1.91(--13) & $\leq$ 7.78(--13) & 160 & $\leq$ 76 & $\leq$ 2.98(--13) & $\leq$ 7.61(--13) \\
& 36 &\ \ \ 3.01(--12) &\ \ \ 1.48(--11) & &\ \ \ \ 85 &\ \ \ 1.61(--9) &\ \ \  1.57(--8) \\
& 37 &\ \ \ 6.93(--11) &\ \ \ 3.55(--10) & &\ \ \ \ 94 &\ \ \ 1.25(--4) &\ \ \  1.17(--3) \\
& 38 &\ \ \ 2.57(--9) &\ \ \ 1.30(--8) & &\ \ \ 103 &\ \ \ 2.64(--3) &\ \ \ 1.51(--1) \\
& 39 &\ \ \ 1.93(--7) &\ \ \ 9.58(--7) & &\ \ \ 112 &\ \ \ 2.35(--3) &\ \ \ 1.16(0) \\
80 & $\leq$ 53 & $\leq$ 2.04(--13) & $\leq$ 6.92(--13) & 320 & $\leq$ 106 & $\leq$ 8.65(--13) & $\leq$ 7.39(--13) \\
& 57 &\ \ \ 2.04(--10) &\ \ \ 5.13(--10) & &\ \ \ 117 &\ \ \ 3.96(--10) &\ \ \ 7.73(--10) \\
& 61 &\ \ \ 3.84(--7) &\ \ \ 9.35(--7) & &\ \ \ 128 &\ \ \ 2.46(--6) &\ \ \ 4.67(--6) \\
& 65 &\ \ \ 1.94(--3) &\ \ \ 4.61(--3) & &\ \ \ 139 &\ \ \ 2.94(--2) &\ \ \ 6.27(--2) \\
& 69 &\ \ \ 1.87(--1) &\ \ \ 6.14(0) & &\ \ \ 150 &\ \ \ 1.15(--3) &\ \ \ 2.18(--2)
\end{tabular}
\medskip

\begin{center}
TABLE 4.1.  Errors in the recursion coefficients $\alpha_k$, $\beta_k$ of (4.5) \\
computed by Stieltjes's procedure
\end{center}

\bigskip

4.3. {\it Multiple-component discretization procedure}.
We assume now a measure $d \lambda$ of the form
$$
d \lambda (t) = w(t) dt + \sum_{j=1}^p y_j \delta (t - x_j ) dt, ~~~~
p \geq 0 ,
\eqno(4.6)
$$
consisting of a continuous part, $w(t) dt$, and (if $p > 0$)
a discrete part written in terms of the Dirac $\delta$-function.
The support of the continuous part is assumed to be an interval or a
finite union of disjoint intervals, some of which may extend to infinity.
In the discrete part, the abscissae $x_j$ are assumed pairwise distinct,
and the weights positive, $y_j > 0$.
The inner product (1.1), therefore, has the form
$$
(u,v) = \int_{\mbox{\myfont R}} u(t) v(t) w(t) dt +
\sum_{j=1}^p y_j u(x_j ) v (x_j ) .
\eqno(4.7)
$$

The basic idea of the discretization procedure
is rather simple:
One approximates the continuous part of the inner product, i.e.,
the integral in (4.7), by a sum, using a suitable quadrature scheme.
If the latter involves $N$ terms, this replaces the inner product (4.7)
by a discrete inner product $( \cdot ~, \cdot )_{N+p}$ consisting of
$N + p$ terms, the $N$ ``quadrature terms'' and the $p$ original terms.
In effect, the measure $d \lambda$ in (4.6) is approximated by a discrete
$(N + p)$-point measure $d \lambda_{N+p}$.
We then compute the desired recursion coefficients from the formula (4.1),
in which the inner product $( \cdot ~, \cdot )$ is replaced, throughout,
by $( \cdot~ , \cdot )_{N+p}$.
Thus, in effect, we approximate
$$
\alpha_k ( d \lambda ) \approx \alpha_k
(d \lambda_{N+p} ) , ~~~~
\beta_k ( d \lambda ) \approx \beta_k ( d \lambda_{N+p} ) .
\eqno(4.8)
$$
The quantities on the right can be computed by the methods of \S4.1 or
\S4.2, i.e., employing the routines {\tt sti} or {\tt lancz}.

The difficult part of this approach is to find a discretization that
results in rapid convergence, as $N \rightarrow \infty$, of the
approximations on the right of (4.8) to the exact values on the left,
even in cases where the weight function $w$ in (4.6) exhibits singular
behavior.
(The speed of convergence, of course, is unaffected
by the discrete part of the inner product (4.7).)
To be successful in this endeavor often requires considerable
inventiveness on the part of the user.
Our routines, {\tt mcdis} and {\tt dmcdis}, that implement this idea
in single resp. double precision, however, are designed to be flexible
enough to promote the use of effective discretization procedures.

Indeed, if the support of the weight function $w$ in (4.7) is contained
in the (finite or infinite) interval ($a,b$), it will often be
useful to first decompose that interval into a finite number of
subintervals,
$$
{\rm supp} ~w \subset [a,b] = \bigcup_{i=1}^m
[a_i , b_i ] , ~~ m \geq 1 ,
\eqno(4.9)
$$
and to approximate the inner product separately on each
subinterval [$a_i , b_i$], using an appropriate weighted
quadrature rule.
Thus, we write the integral in (4.7) as
$$
\int_{\mbox{\myfont R}} u(t) v(t) w(t) dt =
\sum_{i=1}^m \int_{a_i}^{b_i} u(t) v(t) w_i (t) dt ,
\eqno(4.10)
$$
where $w_i$ is an appropriate weight function on [$a_i , b_i$].
The intervals [$a_i , b_i$]
are not necessarily disjoint.
For example, the weight function $w$ may be the sum
$w = w_1 + w_2$ of two weight functions on $[a,b]$, which we
may want to treat individually (cf. Example 4.2 below).
In that case, one would take
$[a_1 , b_1 ] = [a_2 , b_2 ] = [a,b]$ and $w_1$ on the first interval,
$w_2$ on the other.
Alternatively, we may simply want to use a composite quadrature rule to
approximate the integral, in which case
(4.9) is a partition of [$a,b$] and
$w_i (t) = w (t)$ for each $i$.
Still another example is a weight function $w$ which is already supported
on a union of disjoint intervals; in this case, (4.9) would be the same
union, or possibly a refined union where some of the
subintervals are further partitioned.

In whichever way (4.9) and (4.10) are constructed, each integral
on the right of
(4.10) is now approximated by an appropriate quadrature rule,
$$
\int_{a_i}^{b_i} u(t) v(t) w_i (t) dt \approx Q_i (uv) ,
\eqno(4.11)
$$
where
$$
Q_i f = \sum_{r=1}^{N_i} w_{r,i} f (x_{r,i} ) .
\eqno(4.12)
$$
This gives rise to the approximate inner product
\renewcommand{\arraystretch}{2}
$$
\begin{array}{c}
{\displaystyle
(u,v)_{N+p} = \sum_{i=1}^ m ~
\sum_{r=1}^{N_i} w_{r,i} u( x_{r,i} ) v( x_{r,i} ) +
\sum_{j=1}^p y_j u (x_j ) v(x_j ) , } \\
{\displaystyle
N = \sum_{i=1}^m N_i . }
\end{array}
\eqno(4.13)
$$
\renewcommand{\arraystretch}{1}In our routine {\tt mcdis}, we
have chosen, for simplicity, all $N_i$ to be the same integer
$N_0$,
$$
N_i = N_0 , ~~~~
i = 1, 2,\ldots, m ,
\eqno(4.14)
$$
so that $N = m N_0$.
Furthermore, if $n$ is the number of $\alpha_k$ and the number
of $\beta_k$ desired, we have used the following iterative procedure
to determine the coefficients $\alpha_k$, $\beta_k$ to a prescribed (relative) accuracy
$\varepsilon$:
We let $N_0$ increase through a sequence
$\{ N_0^{[s]} \}_{s = 0,1,2,\ldots}$ of integers,
for each $s$ use Stieltjes's (or Lanczos's) algorithm to
compute $\alpha_k^{[s]} = \alpha_k ( d \lambda_{m N_{0}^{[s]} + p} )$,
$\beta_k^{[s]} = \beta_k ( d \lambda_{m N_{0}^{[s]} + p} )$,
$k = 0, 1,\ldots, n-1$, and stop the iteration for the first
$s \geq 1$ for which all inequalities
$$
| \beta_k^{[s]} - \beta_k^{[s-1]} | \leq
\varepsilon \beta_k^{[s]} , ~~~~
k = 0, 1,\ldots, n-1 ,
\eqno(4.15)
$$
are satisfied.
An error flag is provided if within a preset range
$N_0^{[s]} \leq N_0^{\max}$ the stopping criterion (4.15) cannot be
satisfied.
Note that the latter is based solely on the $\beta$-coefficients.
This is because, unlike the $\alpha$'s, they are known to be always
positive, so that it makes sense to insist on relative accuracy.
(In our routine we actually replaced $\beta_k^{[s]}$ on the right
of (4.15) by its absolute value to insure proper termination in cases
of sign-changing measures $d \lambda$.)

In view of the formulae (4.1), it is reasonable to expect, and indeed
has been observed in practice, that satisfaction of (4.15) entails
sufficient absolute accuracy for the $\alpha$'s if they are zero or small,
and relative accuracy otherwise.

Through a bit of experimentation, we have settled on the following
sequence of integers $N_0^{[s]}$:
\renewcommand{\arraystretch}{2}
$$
\begin{array}{c}
{\displaystyle
N_0^{[0]} = 2n , ~~
N_0^{[s]} = N_0^{[s-1]} + \Delta_s , ~~ s = 1, 2,\ldots, } \\
{\displaystyle
\Delta_1 = 1, ~~
\Delta_s = 2^{\lfloor s/5 \rfloor} \cdot n,~~ s = 2, 3,\ldots ~.}
\end{array}
\eqno(4.16)
$$
\renewcommand{\arraystretch}{1}

It should be noted that if the quadrature formula (4.11) is exact
for each $i$, whenever $u \cdot v$ is a polynomial of degree
$\leq 2n-1$ (which is the maximum degree occurring in the inner products
of (4.1), when $k \leq n-1$), then our procedure converges
after the very first iteration step!
Therefore, if each quadrature rule $Q_i$ has (algebraic) degree of
exactness $\geq d (N_0 )$, and
$d(N_0 ) / N_0 = \delta + O ( N_0^{-1} )$ as
$N_0 \rightarrow \infty$, then
we let $N_0^{[0]} = 1 + \lfloor (2n-1 ) / \delta \rfloor $ in an attempt
to get exact answers after one iteration.
Normally, $\delta = 1$ (for interpolatory rules) or $\delta = 2$
(for Gauss-type rules).

The calling sequence of the multiple-component discretization routine is

\medskip
\noindent
\hspace{.75in}{\tt mcdis(n,ncapm,mc,mp,xp,yp,quad,eps,iq,idelta,irout,}

\noindent
\hspace{.75in}{\tt \ \ finl,finr,endl,endr,xfer,wfer,alpha,beta,ncap,}

\noindent
\hspace{.75in}{\tt \ \ kount,ierr,ie,be,x,w,xm,wm,p0,p1,p2)}

\noindent
\hspace{.75in}{\tt dimension xp(*),yp(*),endl(mc),endr(mc),xfer(ncapm),}

\noindent
\hspace{.75in}{\tt \ \ wfer(ncapm),alpha(n),beta(n),be(n),x(ncapm),}

\noindent
\hspace{.75in}{\tt \ \ w(ncapm),xm(*),wm(*),p0(*),p1(*),p2(*)}

\noindent
\hspace{.75in}{\tt logical finl,finr}
\medskip

\noindent
\begin{tabular}{lp{3in}}
On entry,  & \\
\\
\ \ \ \ \ {\tt n} & is the number of recursion coefficients
desired; type integer \\
\\
\ \ \ \ \ {\tt ncapm} & is the integer $N_0^{\max}$ above, i.e., the
maximum integer $N_0$ allowed
({\tt ncapm} = 500 will usually be satisfactory) \\
\\
\ \ \ \ \ {\tt mc} & is the number of component intervals
in the continuous part of the spectrum; type integer \\
\\
\ \ \ \ \ {\tt mp} & is the number of points in the discrete part
of the spectrum; type integer;
if the measure has no discrete part, set {\tt mp} = 0 \\
\\
\ \ \ \ \ {\tt xp,yp} & are arrays of dimension {\tt mp} containing the
abscissae and the jumps of the point spectrum \\
\\
\ \ \ \ \ {\tt quad} & is a subroutine determining the discretization
of the inner product on each component interval, or a dummy routine if
{\tt iq} $\neq$ 1 (see below); specifically, {\tt quad(n,x,w,i,ierr)}
produces the abscissae {\tt x}$(r) = x_{r,i}$ and weights
\end{tabular}

\noindent
\begin{tabular}{lp{3in}}
& \\
\ \ \ \ \hbox to 1in{\ } & 
{\tt w}$(r) = w_{r,i}$, $r = 1, 2,\ldots, n$, of the $n$-point
discretization of the inner product on the interval
$[ a_i , b_i ]$ (cf. (4.13));
an error flag {\tt ierr} is provided to signal the occurrence of an
error condition in the quadrature process \\
\\
\ \ \ \ \ {\tt eps} & is the desired relative accuracy of the nonzero
recursion coefficients; type real \\
\\
\ \ \ \ \ {\tt iq} & is an integer selecting a user-supplied quadrature
routine {\tt quad} if {\tt iq} = 1 or the ORTHPOL routine {\tt qgp}
(see below) otherwise \\
\\ 
\ \ \ \ \ {\tt idelta} & is a nonzero integer, typically 1 or 2,
inducing fast convergence in the case of special quadrature routines;
the default value is {\tt idelta} = 1 \\
\\ 
\ \ \ \ \ {\tt irout} & is an integer selecting the routine for
generating the recursion coefficients from the discrete inner product;
specifically, {\tt irout} = 1 selects the routine {\tt sti}, and
{\tt irout} $\neq$ 1 the routine {\tt lancz}.
\end{tabular}

\bigskip
\noindent
The logical variables {\tt finl,finr} and the arrays {\tt endl,endr,xfer,wfer}
are input variables to the subroutine {\tt qgp} and are
used (and hence need to be properly dimensioned) only if {\tt iq}
$\neq$ 1.

\medskip
\noindent
\begin{tabular}{lp{3in}}
On return, & \\
\\
\ \ \ \ \ {\tt alpha,beta} & are arrays of dimension {\tt n} holding
the desired recursion coefficients {\tt alpha} $(k) = \alpha_{k-1}$,
{\tt beta} $(k) = \beta_{k-1}$, $k = 1, 2,\ldots,$ {\tt n} \\
\\
\ \ \ \ \ {\tt ncap} & is the integer $N_0$ yielding  
convergence \\
\\
\ \ \ \ \ {\tt kount} & is the number of iterations required to
achieve convergence \\
\end{tabular}

\noindent
\begin{tabular}{lp{3in}}
\ \ \ \ \ {\tt ierr} & is an error flag, equal to 0 on normal return,
equal to --1 if {\tt n} is not in the proper range,
equal to $i$ if there is an error condition in the discretization
on the $i$th interval, and equal to {\tt ncapm} if the discretized
Stieltjes procedure does not converge within the discretization
resolution specified by {\tt ncapm} \\
\\
\ \ \ \ \ {\tt ie} & is an error flag inherited from the routine
{\tt sti} or {\tt lancz} (whichever is used).
\end{tabular}
\medskip

\noindent
The arrays {\tt be,x,w,xm,wm,p0,p1,p2} are used for working space,
the last five having dimension {\tt mc}$\times${\tt ncapm}+{\tt mp}.

A general-purpose quadrature routine, {\tt qgp}, is provided for cases
in which it may be difficult to develop special discretizations that take advantage
of the structural properties of the weight function $w$ at hand.
The routine assumes the same setup (4.9)--(4.14) used in
{\tt mcdis}, with {\it disjoint} intervals [$a_i , b_i$],
and provides for $Q_i$ in (4.12) the Fej\'{e}r quadrature
rule, suitably transformed to the interval
[$a_i , b_i$], with the same number
$N_i = N_0$ of points for each $i$.
Recall that the $N$-point Fej\'{e}r rule on the standard interval [--1,1] is the
interpolatory quadrature formula
$$
Q_N^F f = \sum_{r=1}^N w_r^F f (x_r^F ) ,
\eqno(4.17)
$$
where $x_r^F = \cos ((2r-1) \pi /2N)$ are the Chebyshev points.
The weights are all positive and can be computed explicitly in terms
of trigonometric functions (cf., e.g., [15]).
The rule (4.17) is now applied to the integral in (4.11) by
transforming the
interval [--1,1] to [$a_i , b_i$] via some monotone function
$\phi_i$ (a linear function if $[a_i , b_i ]$ is finite)
and letting $f = u v w_i$:
$$
\int_{a_i}^{b_i} u(t) v(t) w_i (t) dt = \int_{-1}^1
u ( \phi_i ( \tau ) ) v ( \phi_i ( \tau ) )
w_i ( \phi_i ( \tau ) ) \phi_i^{\prime} ( \tau ) d \tau
$$
$$
\approx \sum_{r=1}^N
w_r^F w_i ( \phi_i ( x_r^F ) ) \phi_i^{\prime} (x_r^F ) \cdot
u ( \phi_i ( x_r^F ) ) v ( \phi_i ( x_r^F ) ) .
$$
Thus, in effect, we take in (4.13)
$$
x_{r,i} = \phi_i ( x_r^F ) , ~~
w_{r,i} = w_r^F w_i ( \phi_i (x_r^F ) ) \phi_i^{\prime}
(x_r^F ) , ~~
i = 1, 2,\ldots, m .
\eqno(4.18)
$$
If the interval [$a_i , b_i$] is half-infinite,
say of the form [0,$\infty$], we use
$\phi_i (t) = (1+t) / (1 - t)$, and
similarly for intervals of the form
[--$\infty$,$b$] and [$a, \infty$].
If $[a_i , b_i ] = [- \infty , \infty ]$, we use
$\phi_i (t) = t / (1 - t^2 )$.

The routine {\tt qgp} has the following calling sequence:

\medskip
\begin{verbatim}
     subroutine qgp(n,x,w,i,ierr,mc,finl,finr,endl,endr,xfer,wfer)
     dimension x(n),w(n),endl(mc),endr(mc),xfer(*),wfer(*)
     logical finl,finr
\end{verbatim}
\medskip

\noindent
\begin{tabular}{lp{3in}}
On entry,  & \\
\\
\ \ \ \ \ {\tt n} & is the number of terms in the Fej\'{e}r
quadrature rule \\
\\
\ \ \ \ \ {\tt i} & indexes the interval [$a_i , b_i$] for which the
quadrature rule is desired; an interval that extends to
--$\infty$ has to be indexed by 1, one that extends
to +$\infty$ by {\tt mc} \\
\\
\ \ \ \ \ {\tt mc} & is the number of component intervals; type integer \\
\\
\ \ \ \ \ {\tt finl} & is a logical variable to be set .{\tt true}. if
the extreme left interval is finite, and .{\tt false}. otherwise \\
\\
\ \ \ \ \ {\tt finr} & is a logical variable to be set .{\tt true}. if
the extreme right interval is finite, and .{\tt false}. otherwise \\

\\
\ \ \ \ \ {\tt endl} & is an array of dimension {\tt mc} containing the
left endpoints of the component intervals; if the first of these
extends to --$\infty$, {\tt endl}(1) is not being used by the routine\\
\\
\ \ \ \ \ {\tt endr} & is an array of dimension {\tt mc} containing
the right endpoints of the component intervals;
if the last of these extends to +$\infty$, {\tt endr}({\tt mc}) is not
being used by the routine 
\end{tabular}

\noindent
\begin{tabular}{lp{3in}}
\ \ \ \ \ {\tt xfer,wfer} & are working arrays 
holding the standard
Fej\'{e}r nodes and weights, respectively; they are dimensioned
in the routine {\tt mcdis}. \\
On return, & \\
\\
\ \ \ \ \ {\tt x,w} & are arrays of dimension {\tt n} holding
the abscissae and weights (4.18) of the discretized inner product for the
$i$th component interval \\
\\
\ \ \ \ \ {\tt ierr} & has the integer value 0.
\end{tabular}
\bigskip

\noindent
The routine calls on the subroutines {\tt fejer,symtr} and
{\tt tr}, which are appended to the routine {\tt qgp} in \S4 of
the package.
The first generates the Fej\'{e}r quadrature
rule, the others perform variable transformations.
The user has to provide his own function routine
{\tt wf(x,i)} to calculate the weight function
$w_i (x)$ on the $i$th component interval.

\bigskip

{\it Example} 4.2.
Chebyshev weight plus a constant:
$w^c (t) = (1-t^2 )^{- 1/2} + c$,
$c > 0$, $-1 < t < 1$.
\smallskip

It would be difficult, here, to find a single quadrature rule for
discretizing the inner product and obtain fast convergence.
However, using in (4.9) $m = 2$,
$[a_1 , b_1 ]$ = $[a_2 , b_2 ]$ = $[-1,1]$,
and $w_1 (t) = (1-t^2)^{- 1/2}$, $w_2 (t) = c$ in (4.11),
and taking for $Q_1$ the Gauss-Chebyshev, and for $Q_2$ the 
Gauss-Legendre $n$-point rule (the latter multiplied by $c$), 
yields convergence to $\alpha_k (w^c )$, 
$\beta_k (w^c )$, $k = 0, 1,\ldots, n-1$, in one iteration
(provided $\delta$ is set equal to 2)! 
Actually, we need $N_0 = n+1$, in order to test for convergence; cf.
(4.15).
The driver {\tt test4} implements this technique and calculates the first $n = 80$
beta-coefficients to a relative accuracy of 5000$\times \varepsilon^{s}$
for $c = 1$, 10, 100.
(All $\alpha_k$ are zero.)
Attached to the driver is the quadrature routine {\tt qchle} used in this
example.
It, in turn, calls for the Gauss quadrature routine {\tt gauss}, to be
described later in \S6.
Anticipating convergence after one iteration, we put
{\tt ncapm} = 81.

The weight function of Example 4.2 provides a continuous link between
the Chebyshev polynomials ($c = 0$) and the Legendre
polynomials ($c = \infty$); the
recursion coefficients $\beta_k (w^c )$ indeed converge (except for $k = 0$)
to those of the Legendre polynomials, as $c \rightarrow \infty$.

Selected results of {\tt test4} (where {\tt irout} in {\tt mcdis}
can be arbitrary) are shown in Table 4.2.
The output variable {\tt kount} is 1 in each
case, confirming convergence after one iteration.
The coefficients $\beta_0 (w^c )$ are easily seen to be
$\pi + 2c$.

\bigskip
\begin{center}
\begin{tabular}{clll}
\multicolumn{1}{c}{$k$} &
\multicolumn{1}{c}{$\beta_k (w^1 )$} &
\multicolumn{1}{c}{$\beta_k (w^{10})$} &
\multicolumn{1}{c}{$\beta_k (w^{100})$} \\ \hline
 & & &  \\*[-10pt]
0 & 5.141592654 & 23.14159265 & 203.1415927 \\
1 & \ \ .4351692451 & \ \ \ .3559592080 & \ \ \ \ \ .3359108398 \\
5 & \ \ .2510395775 & \ \ \ .2535184776 & \ \ \ \ \ .2528129500 \\
12 & \ \  .2500610870 & \ \ \ .2504824840 & \ \ \ \ \  .2505324193 \\
25 & \ \ .2500060034 & \ \ \ .2500682357 & \ \ \ \ \ .2501336338 \\
51 & \ \ .2500006590 & \ \ \ .2500082010 & \ \ \ \ \ .2500326887 \\
79 & \ \ .2500001724 & \ \ \ .2500021136 & \ \ \ \ \ .2500127264
\end{tabular}
\end{center}
\medskip
\begin{center}
TABLE 4.2.  Selected recursion coefficients $\beta_k ( w^c )$ for $c$ = 1, 10, 100
\end{center}
\bigskip

{\it Example} 4.3.
Jacobi weight with one mass point at the left endpoint:
$w^{( \alpha , \beta )} (x~;y)$ = $[ \mu_0^{( \alpha , \beta )}]^{-1}
(1-x)^{\alpha} (1+x)^{\beta}$ + $y~ \delta (x+1)$ on ($-$1,1),
$\mu_0^{( \alpha , \beta )} = 2^{ \alpha + \beta +1}
\Gamma ( \alpha + 1) \Gamma ( \beta +1 ) / \Gamma ( \alpha + \beta + 2 )$,
$\alpha > - 1$, $\beta > - 1$, $y > 0$.
\smallskip

The recursion coefficients $\alpha_k$, $\beta_k$ are known explicitly
(see [6, Eqs. (6.23), (3.5)]\footnote{
In [6], the interval is taken to be [0,2], rather than [$-$1,1].
There is a typographical error in the first formula of (6.23), which
should have the numerator $2 \beta + 2$ instead of $2 \beta +1$.})
and can be expressed, with some effort, in terms of the
recursion coefficients $\alpha_k^J$, $\beta_k^J$ for the Jacobi
weight
$w^{( \alpha , \beta )} ( \cdot ) = w^{( \alpha , \beta )}$
$( \cdot ~; 0)$.
The formulae are:

\renewcommand{\arraystretch}{2}
$$
\begin{array}{c}
{\displaystyle
\alpha_0 = ~\frac{\alpha_0^J - y}{1 + y} , ~~~~
\beta_0 = \beta_0^J + y , } \\
\alpha_k = \alpha_k^J + ~\frac{2k( \alpha + k )}
{( \alpha + \beta + 2k)( \alpha + \beta + 2 k+1)} ~
(c_k - 1) + ~\frac{2( \beta +k+1)( \alpha + \beta +k+1)}
{( \alpha + \beta +2k+1)( \alpha + \beta + 2k+2)} ~
\left( \frac{1}{c_k} - 1 \right) , \\
{\displaystyle
\beta_k = ~\frac{c_k}{c_{k-1}} ~\beta_k^J , ~~~~
k = 1, 2, 3,\ldots, }
\end{array}
\eqno(4.19)
$$
\renewcommand{\arraystretch}{1}where
$$
c_0 = 1 + y , ~~
c_k = ~
\frac{ 1 + ~\frac{( \beta +k+1)( \alpha + \beta + k+1)}{k ( \alpha + k )}
y d_k }
{1 + yd_k } ~, ~~
k = 1, 2,\ldots,
\eqno(4.20)
$$
and
\renewcommand{\arraystretch}{2}
$$
\begin{array}{l}
{\displaystyle
d_1 = 1,} \\
{\displaystyle
d_k = \frac{( \beta +k)( \alpha + \beta +k)}
{( \alpha + k-1)(k-1)} ~ d_{k-1} , ~~~~
k = 2, 3,\ldots ~. }
\end{array}
\eqno(4.21)
$$
\renewcommand{\arraystretch}{1}

Again, it is straightforward with {\tt mcdis} to get exact results
(modulo rounding) after one iteration, by using the Gauss-Jacobi
quadrature rule (see {\tt gauss} in \S6) to discretize the continuous
part of the measure.
The driver {\tt test5} generates in this manner the first
$n$ = 40 recursion coefficients
$\alpha_k$, $\beta_k$, $k = 0, 1,\ldots, n-1$,
to a relative accuracy of 5000$\times \varepsilon^{s}$,
for $y = \frac{1}{2}$, 1, 2, 4, 8.
For each $\alpha = -.8(.2)1.$ and
$\beta$ = --.8(.2)1., it computes the maximum relative errors
(absolute error, if $\alpha_k \approx 0$) of the $\alpha_k$,
$\beta_k$ by comparing them with the exact coefficients.
These have been computed in double precision
by a straightforward implementation of the formulae
(4.19)--(4.21). 

As expected, the output of {\tt test5} reveals convergence after one
iteration, the variable {\tt kount}
having consistently the value 1.
The maximum relative error in the $\alpha_k$ is found to generally lie
between $2 \times 10^{-8}$ and $3 \times 10^{-8}$, the one in the
$\beta_k$ between $7.5 \times 10^{-12}$ and $8 \times 10^{-12}$;
they are attained for $k$ at or near 39.
The discrepancy between the errors in the $\alpha_k$ and those in the
$\beta_k$ is due to the $\alpha_k$ being considerably smaller than the
$\beta_k$, by 3--4 orders of magnitude.
Replacing the routine {\tt sti} in {\tt mcdis} by
{\tt lancz} yields very much the same error picture.

It is interesting to note that the addition of a second mass point at
the other endpoint makes an analytic determination of the recursion
coefficients intractable (cf. [6, p. 713]).
Numerically, however, it makes no difference whether there are two or more
mass points, and whether they are located inside, or outside, or on the
boundary of the support interval.
It was observed, however, that if at least one mass point is
located outside the interval [$-$1,1], the
procedure {\tt sti} used in
{\tt mcdis} becomes severely
unstable\footnote{
This has also been observed in the similar Example 4.8 of [19], but was
incorrectly attributed to a phenomenon of ill-conditioning.
Indeed, the statement made at the end of Example 4.8 can now be retracted:
stable methods {\it do} exist, namely the method embodied by the routine
{\tt mcdis} in combination with {\tt lancz}.}
and {\it must} be replaced by {\tt lancz}.

\bigskip

{\it Example} 4.4.
Logistic density function: $w(x) = e^{-x} / (1 + e^{-x} )^2$ on
($- \infty, \infty$).
\smallskip

In this example we illustrate a slight variation of the discretization
procedure (4.9)--(4.13), which ends up with a discrete inner product of the
same type as in (4.13) (and thus implementable by the routine
{\tt mcdis}) but derived in a somewhat different manner.
The idea is to integrate functions with respect to the density $w$
by splitting the integral into two parts, one from $- \infty$
to 0, the other from 0 to $\infty$, changing variables in the first part,
and thus obtaining
$$
\int_{- \infty}^\infty f(t) w(t) dt =
\int_0^\infty f(-t) ~
\frac{e^{-t}}{(1 + e^{-t} )^2} ~
dt + \int_0^\infty f(t) ~
\frac{e^{-t}}{(1 + e^{-t} )^2 } ~dt .
\eqno(4.22)
$$
Since $(1 + e^{-t} )^{-2}$ quickly tends to 1 as $t \rightarrow \infty$,
a natural
discretization of both integrals is provided by the Gauss-Laguerre
quadrature rule applied to the product $f( \pm t ) \cdot$
$(1 + e^{-t} )^{-2}$.
This amounts to taking, in (4.13), $m = 2$ and
$$
x_{r,1} = - x_r^L , ~
x_{r,2} = x_r^L ; ~~
w_{r,1} = w_{r,2} = ~\frac{w_r^L}{(1 + e^{- x_r^L} )^2} ~, ~
r = 1, 2,\ldots, N ,
$$
where $x_r^L$, $w_r^L$ are the Gauss-Laguerre $N$-point
quadrature nodes and weights.

The driver {\tt test6} incorporates this discretization into the routines
{\tt mcdis} and {\tt dmcdis}, runs them for $n = 40$ with error
tolerances 5000$\times \varepsilon^{s}$
and 1000$\times \varepsilon^{d}$, respectively,
and prints the absolute errors in the $\alpha$'s
($\alpha_k = 0$, in theory) and the relative errors in the $\beta$'s.
(We used the default value $\delta = 1$.)
Also printed are the number of iterations \#it (= {\tt kount}) in (4.15)
and the corresponding final value $N_0^f$ (= {\tt ncap}).
In single precision,
we found \#it = 1, $N_0^f = 81$, and in double precision,
\#it = 5, $N_0^f = 281$.
Both routines returned with the error flags equal to 0, indicating a 
normal course 

\bigskip
\begin{center}
\begin{tabular}{rlll}
\multicolumn{1}{c}{$k$} &
\multicolumn{1}{c}{$\beta_k$} &
\multicolumn{1}{c}{err $\alpha_k$} &
\multicolumn{1}{c}{err $\beta_k$} \\ \hline
 & & &  \\*[-10pt]
0  & 1.000000000000000000000000 & 4.572(--13) & 1.918(--13) \\
1  & 3.289868133696452872944830 & 1.682(--13) & 5.641(--13) \\
6  & 89.44760352315950188817832 & 2.187(--12) & 2.190(--12) \\
15 & 555.7827839879296775066697 & 1.732(--13) & 2.915(--12) \\
26 & 1668.580222268668421827788 & 3.772(--12) & 4.112(--12) \\
39 & 3753.534025194898387722354 & 2.482(--11) & 4.533(--12)
\end{tabular}
\end{center}
\smallskip
\begin{center}
TABLE 4.3.  Selected output from {\tt test6}
\end{center}
\medskip

\noindent
of events.
A few selected double-precision values of the coefficients $\beta_k$
along with absolute errors in the $\alpha$'s and relative errors in
the  $\beta$'s are shown in Table 4.3. The results are essentially
the same no matter whether {\tt sti} or {\tt lancz} is used in
{\tt mcdis}.
The maximum errors observed are
$2.482 \times 10^{-11}$ for the $\alpha$'s, and $4.939 \times 10^{-12}$
for the $\beta$'s, which are well within the single-precision tolerance
$\varepsilon = 5000 \times \varepsilon^{s}$.

On computers with limited exponent range, convergence difficulties
may arise, both with {\tt sti} and {\tt lancz},
owing to underflow in many of the Laguerre quadrature
weights. This seems to perturb the problem significantly enough to
prevent the discretization procedure from converging.
\bigskip

{\it Example} 4.5.
Half-range Hermite measure: $w(x) = e^{-x^{2}}$ on ($0, \infty$).
\smallskip

This is an example of a measure for which there do not seem to exist natural
discretizations other than those based on
composite quadrature rules.
Therefore, we applied our general-purpose routine
{\tt qgp} (and its double-precision companion {\tt dqgp}),
using, after some experimentation, the partition
$[0, \infty ] = [0,3] \cup [3,6] \cup [6,9] \cup [9, \infty ]$.
The driver {\tt test7} implements this, with $n = 40$ and an error 
tolerance 50$\times \varepsilon^{s}$ in single precision,
and 1000$\times \varepsilon^{d}$ in double precision.

The single-precision routine {\tt mcdis} (using the default
value $\delta = 1$) converged after one
iteration, returning {\tt ncap} = 81, whereas the double-precision

\begin{center}
\begin{tabular}{rll}
\multicolumn{1}{c}{$k$} &
\multicolumn{1}{c}{$\begin{array}{c} \alpha_k \\ \mbox{err}~ \alpha_k \end{array}$} &
\multicolumn{1}{c}{$\begin{array}{c} \beta_k \\ \mbox{err}~ \beta_k \end{array}$} \\ \hline
 & &  \\*[-10pt]
0 & \ \ .5641895835477562869480795 & \ \ .8862269254527580136490837 \\
& 1.096(--13) & 3.180(--13) \\
1 & \ \ .9884253928468002854870634 & \ \ .1816901138162093284622325 \\
& 1.514(--13) & 7.741(--14) \\
6 & 2.080620336400833224817622 & 1.002347851011010842224538 \\
& 1.328(--13) & 5.801(--14) \\
15 & 3.214270636071128227448914 & 2.500927917133702669954321 \\
& 2.402(--14) & 8.186(--14) \\
26 & 4.203048578872001952660277 & 4.333867901229950443604430 \\
& 1.415(--13) & 7.878(--14) \\
39 & 5.131532886894296519319692 & 6.500356237707132938035155 \\
& 6.712(--13) & 1.820(--14)
\end{tabular}
\end{center}
\smallskip

\begin{center}
TABLE 4.4.  Selected output from {\tt test7}
\end{center}
\bigskip

\noindent
routine {\tt dmcdis} took four iterations
to converge, and returned {\tt ncapd} = 201.
Selected results (where err $\alpha_k$ and err $\beta_k$
both denote relative errors) are shown in Table 4.4. The maximum error
err $\alpha_k$ occurred at $k = 10$, and had the
value $1.038 \times 10^{-12}$, whereas $\max_k$ err $\beta_k$ =
$3.180 \times 10^{-13}$ is
attained at $k = 0$.
The latter is within the error tolerance $\varepsilon$, the
former only slightly larger.
Comparison of the double-precision results with Table I on the microfiche
supplement to [12] revealed agreement to all 20 decimal digits given there,
for all $k$ in the range $0 \leq k \leq 19$. Interestingly, the
routine {\tt sti} in {\tt mcdis} did consistently better than
{\tt lancz} on the $\beta$'s, by a factor as large as 235 (for
$k = 33$), and is comparable with {\tt lancz} (sometimes better,
sometimes worse) on the $\alpha$'s.

Without composition, i.e.,
using {\tt mc}=1 in {\tt mcdis}, it takes 8 iterations
$(N_0^f = 521)$ in single precision, and 10 iterations
$(N_0^f = 761)$ in double precision, to satisfy the much weaker error
tolerances $\varepsilon = \frac{1}{2}~10^{-6}$ and
$\varepsilon^{d} = \frac{1}{2}~10^{-12}$, respectively.
All single-precision results, however, turn out to be accurate to about
12 decimal places.
(This is because of the relatively large final increment
$\Delta_8 = 2n = 80$ in $N_0$ (cf. (4.16)) that forces convergence.)
\bigskip

4.4.  {\it Discretized modified Chebyshev algorithm}.
The whole apparatus of discretization (cf. (4.9)--(4.14)) can also
be employed in connection with the modified Chebyshev algorithm
(cf. \S3), if one discretizes modified moments rather than inner
products.
Thus, one approximates (cf. (4.14), (4.16))
$$
\nu_k ( d \lambda ) \approx \nu_k ( d \lambda_{mN_{0}^{[s]} + p} )
\eqno(4.23)
$$
and iterates the modified Chebyshev algorithm with
$s = 0, 1, 2,\ldots$ until the convergence criterion (4.15) is
satisfied.
(It would be unwise to test convergence on the modified moments,
for reasons explained in [19, \S2.5].)
This is implemented in the routine {\tt mccheb}, whose calling
sequence is
\medskip

\begin{verbatim}
     mccheb(n,ncapm,mc,mp,xp,yp,quad,eps,iq,idelta,finl,
       finr,endl,endr,xfer,wfer,a,b,fnu,alpha,beta,ncap,
       kount,ierr,be,x,w,xm,wm,s,s0,s1,s2)
\end{verbatim}
\medskip

\noindent
Its input and output parameters have the same meaning as in the
routine {\tt mcdis}.
In addition, the arrays {\tt a,b} of dimension 2$\times${\tt n}--1
are to be supplied with the recursion coefficients
{\tt a}$(k) = a_{k-1}$, {\tt b}$(k) = b_{k-1}$,
$k = 1, 2,\ldots,$ 2$\times${\tt n}--1, defining the modified moments.
The arrays {\tt be,x,w,xm,wm,s,s0,s1,s2} are used for working
space.
The double-precision version of the routine has the name {\tt
dmcheb}.

The discretized modified Chebyshev algorithm must be expected to
behave similarly as its close relative, the modified Chebyshev
algorithm.
In particular, if the latter suffers from ill-conditioning, so does
the former.
\bigskip

{\it Example} 4.6.  Example 3.1, revisited.
\smallskip

We recompute the $n = 40$ first recursion coefficients
$\alpha_k$, $\beta_k$ of Example 3.1 to an accuracy of
100$\times \varepsilon^{s}$ in single precision, using the
routine {\tt mccheb} instead of the routine {\tt cheb}.
For the discretization of the modified moments we employed the
Gauss-Chebyshev quadrature rule:
$$
\int_{-1}^1 f(t) (1 - \omega^2 t^2 )^{- 1/2} (1-t^2 )^{- 1/ 2}
dt \approx ~\frac{\pi}{N} ~
\sum_{r=1}^N f (x_r ) (1 - \omega^2 x_r^2 )^{- 1/2} ,
\eqno(4.24)
$$
where $x_r = \cos ((2r-1) \pi /2N )$ are the Chebyshev points.
This is accomplished by the driver {\tt test8}.
The results of this test (shown in the package) agree to all 10 decimal 
places with those of {\tt test1}.
The routine {\tt mccheb} converged in one iteration, with {\tt
ncap} = 81, for $\omega^2 = .1$, .3, .5, .7, .9, in 4 iterations,
with
{\tt ncap} = 201, for $\omega^2 = .99$, and in 8 iterations, with 
{\tt ncap} = 521, for $\omega^2 = .999$.
A double-precision version of {\tt test8} was also run with
$\varepsilon ~ = ~\frac{1}{2} ~ \times 10^{-20}$ (not shown in the package)
and produced correct results to 20 decimals in one iteration
({\tt ncap} = 81) for $\omega^2 = .1$, .3,
.5, .7, in 3 iterations ({\tt ncap} = 161) for 
$\omega^2 = .9$, in 6 iterations ({\tt ncap} = 361) for
$\omega^2 = .99$, and in 11 iterations ({\tt ncap} = 921) for
$\omega^2 = .999$.

\vspace{.375in}
\noindent
{\small 5.  MODIFICATION ALGORITHMS}

\bigskip

Given a positive measure $d \lambda (t)$ supported on the real line,
and two polynomials
$u(t) = \pm~ \Pi_{\rho =1}^r (t - u_{\rho} )$,
$v(t) = \Pi_{\sigma =1}^s (t - v_\sigma )$ whose ratio is finite
on the support of $d \lambda$, we may ask for the
recursion coefficients
$\hat{\alpha}_k = \alpha_k (d \hat{\lambda} )$,
$\hat{\beta}_k = \beta_k ( d \hat{\lambda} )$ of the modified measure
$$
d \hat{\lambda} (t) = \frac{u(t)}{v(t)} ~
d \lambda (t) , ~~~~
t \in \mbox{supp} ( d \lambda ) ,
\eqno(5.1)
$$
assuming known the recursion coefficients
$\alpha_k = \alpha_k ( d \lambda )$,
$\beta_k = \beta_k ( d \lambda )$ of the given measure.
Methods that accomplish the passage from the $\alpha$'s and $\beta$'s
to the $\hat{\alpha}$'s and $\hat{\beta}$'s will be called
{\it modification algorithms}.
The simplest case $s = 0$ (i.e., $v(t) \equiv 1$) and $u$ positive
on supp$(d \lambda )$ has already been considered
by Christoffel [7], who represented the polynomial
$u ( \cdot ) \hat{\pi}_k ( \cdot )$ = $u ( \cdot ) \pi_k$
$( \cdot ~; d \hat{\lambda} )$ in determinantal form in terms of the
polynomials
$\pi_j ( \cdot ) = \pi_j ( \cdot ~; d \lambda )$,
$j = k , k+1 ,\ldots, k+r$.
This is now known as {\it Christoffel's theorem}.
Christoffel, however, did not address the problem of how to generate
the new coefficients
$\hat{\alpha}_k$, $\hat{\beta}_k$ in terms of the old ones.
For the more general modification (5.1), Christoffel's theorem has
been generalized by Uvarov [48,49].
The coefficient problem stated above, in this general case, has been
treated in [20], and previously by Galant [13] in the
special case $v(t) \equiv 1$.

The passage from $d \lambda$ to $d \hat{\lambda}$ can be carried out in
a sequence of elementary steps involving real linear factors
$ t-x $, or real quadratic factors
$(t - x)^2 + y^2$, either in $u(t)$ or in $v(t)$.
The corresponding elementary steps in the passage
from the $\alpha$'s and $\beta$'s to the $\hat{\alpha}$'s and $\hat{\beta}$'s
can all be performed by means of certain nonlinear recurrences.
Some of these, however, when divisions of the measure $d \lambda$
are involved, are liable to instabilities.
An alternative method can then be used, which appeals to the
modified Chebyshev algorithm supplied with appropriate
modified moments.
These latter are of independent interest and find application, e.g.,
in evaluating the kernel in the contour integral representation of
the Gauss quadrature remainder term.
\bigskip

5.1 {\it Nonlinear recurrence algorithms}.
The routine that carries out the elementary modification steps is called
{\tt chri} and has the calling sequence
\medskip

\begin{center}
{\tt chri(n,iopt,a,b,x,y,hr,hi,alpha,beta,ierr).}
\end{center}
\medskip

\noindent
\begin{tabular}{lp{3in}}
On entry,  & \\
\\
\ \ \ \ \ {\tt n} & is the number of recursion coefficients desired;
type integer 
\end{tabular}

\noindent
\begin{tabular}{lp{3in}}
\\
\ \ \ \ \ {\tt iopt} & is an integer identifying the type of modification
as follows: \\
& 1: \ \ $d \hat{\lambda} (t) = ( t-x ) d \lambda (t)$ \\
& 2: \ \ $d \hat{\lambda} (t) = ((t-x)^2 + y^2 ) d \lambda (t)$, $y > 0$  \\
& 3: \ \ $d \hat{\lambda} = (t^2 + y^2 ) d \lambda (t)$ with $d \lambda (t)$ and \\
& \hspace{.25in} supp$(d \lambda )$ assumed symmetric with \\
& \hspace{.25in} respect to the origin and $y > 0$ \\
& 4: \ \ $d \hat{\lambda} (t) = d \lambda (t) / ( t-x )$ \\
& 5: \ \ $d \hat{\lambda} (t) = d \lambda (t) / ((t-x)^2 + y^2 )$, $y > 0$ \\
& 6: \ \ $d \hat{\lambda} (t) = d \lambda (t) / (t^2 + y^2 )$ with $d \lambda (t) $ and \\
& \hspace{.25in} supp$(d \lambda )$ assumed symmetric with \\
& \hspace{.25in} respect to the origin and $y > 0$ \\
& 7: \ \ $d \hat{\lambda} (t) = (t-x)^2 d \lambda (t)$ \\
\ \ \ \ \ {\tt a,b} & are arrays of dimension {\tt n}+1 holding 
the recursion coefficients {\tt a}$(k) = \alpha_{k-1} (d \lambda )$,
{\tt b}$(k) = \beta_{k-1} ( d \lambda )$, $k = 1, 2,\ldots,$ {\tt n}+1
\\
\ \ \ \ \ {\tt x,y} & are real parameters defining the linear
and quadratic factors (or divisors) of $d \lambda$ \\
\\
\ \ \ \ \ {\tt hr,hi} & are the real and imaginary part, respectively,
of $\int_{\mbox{\myfont R}} d \lambda (t) / (z-t)$, where
$z = x+iy$; the parameter {\tt hr} is used only if {\tt iopt} = 4 or
5, the parameter {\tt hi} only if {\tt iopt} = 5 or 6. 
\end{tabular}

\noindent
\begin{tabular}{lp{3in}}
\\
On return, & \\
\\
\ \ \ \ \ {\tt alpha,beta} & are arrays of dimension {\tt n}
containing
the desired recursion coefficients {\tt alpha}$(k) = \alpha_{k-1}$
$(d \hat{\lambda} )$, {\tt beta}$(k) = \beta_{k-1}$ $(d \hat{\lambda} )$,
$k = 1, 2,\ldots,$ {\tt n} \\
\\
\ \ \ \ \ {\tt ierr} & is an error flag, equal to 0 on normal return,
equal to 1 if {\tt n} $\leq 1$ (the routine assumes that {\tt n} is larger
than or equal to 2), and equal to 2 if the integer {\tt iopt} is
inadmissible.
\end{tabular}
\medskip

It should be noted that in the cases {\tt iopt} = 1 and {\tt iopt} = 4,
the modified measure $d \hat{\lambda}$ is positive [negative] definite
if $x$ is to the left [right] of the support of $d \lambda$, but
indefinite otherwise. Nevertheless, it is permissible to have $x$ inside
the support of $d \lambda$ (or inside its convex hull), provided the
resulting measure $d \hat{\lambda}$ is still quasi-definite (cf. [20]).

For {\tt iopt} = $1, 2,\ldots, 6$, the methods used in {\tt chri}
are straightforward implementations of the nonlinear recurrence algorithms respectively
in Eq. (3.7), (4.7), (4.8), (5.1), (5.8) and (5.9) of [20].
The only minor modification required concerns
$\hat{\beta}_0 = \beta_0 ( d \hat{\lambda} )$.
In [20], this constant was taken to be 0, whereas here it is defined
to be $\hat{\beta}_0 = \int_{\mbox{\myfont R}} d \hat{\lambda} (t)$.
Thus, for example, if {\tt iopt} = 2,
$$
\hat{\beta}_0 = \int_{\mbox{\myfont R}} ((t-x)^2 + y^2 )
d \lambda (t) = \int_{\mbox{\myfont R}} (( t - \alpha_0 +
\alpha_0 - x )^2 + y^2  ) d \lambda (t)
$$
$$
= \int_{\mbox{\myfont R}} (( t - \alpha_0 )^2 +
( \alpha_0 - x)^2 + y^2 ) d \lambda (t) ,
$$
since $\int_{\mbox{\myfont R}} (t - \alpha_0 ) d \lambda (t)$ =
$\int_{\mbox{\myfont R}} \pi_1 (t) d \lambda (t) = 0$.
Furthermore (cf. (4.1)),
$$
\int_{\mbox{\myfont R}} (t - \alpha_0 )^2 d \lambda (t)
= \beta_0 \beta_1 ,
$$
so that the formula to be used for
$\hat{\beta}_0$ is
\medskip

\begin{center}
$\hat{\beta}_0 = \beta_0 ( \beta_1 +
( \alpha_0 - x)^2 + y^2 ) ~~~~$
({\tt iopt} = 2) .
\end{center}
\medskip
\noindent
Similar calculations need to be made in the other cases.

The case {\tt iopt} = 7 incorporates a {\it QR} step with
shift $x$, following Kautsky and Golub [42], and uses an
adaptation of the algorithm in [51, Eq. (67.11), p. 567] to carry 
out the {\it QR} step.
The most significant modification made is the replacement of the test
$c \neq 0$ by $| c | > \varepsilon$, where $\varepsilon = 5 \times \varepsilon^{s}$ is a quantity
close to, but slightly larger than, the machine precision.
(Without this modification, the algorithm could fail.)

The methods used in {\tt chri} are believed to be quite
stable when the measure $d \lambda$ is modified multiplicatively
({\tt iopt} = 1, 2, 3 and 7).
When divisions are involved ({\tt iopt} = 4, 5 and 6),
however, the algorithms rapidly become unstable as the point
$z = x+iy \in \mbox{\myfont C}$ moves away from the support
interval of $d \lambda$.
(The reason for this instability is not well understood at present;
see, however, Galant [14].)
For such cases there is an alternative routine, {\tt gchri} (see
\S5.2), that can be used.
\newpage

{\it Example} 5.1.
Checking the results (for $\sigma = \frac{1}{2}$) of {\tt test2}.
\smallskip

We apply {\tt chri} (and the corresponding double precision routine
{\tt dchri}) with {\tt iopt} = 1, $x = 0$, to
$d \lambda_\sigma (t) = t^{\sigma}$ ln $(1/t)$ on (0,1) with
$\sigma = - \frac{1}{2}$, to recompute the results of {\tt test2}
for $\sigma = \frac{1}{2}$.
This can be done by a minor 
\bigskip

\begin{center}
\begin{tabular}{crllll}
\multicolumn{1}{c}{$\sigma$} &
\multicolumn{1}{c}{$k$} &
\multicolumn{1}{c}{err $\alpha_k$} &
\multicolumn{1}{c}{err $\beta_k$} &
\multicolumn{1}{c}{err $\alpha_k^d$} &
\multicolumn{1}{c}{err $\beta_k^d$} \\ \hline
 & & & & &  \\*[-10pt]
.5 & 0 & 7.895(--14) & 4.796(--14) & 2.805(--28) & 7.952(--28) \\
& 12 & 3.280(--12) & 6.195(--12) & 8.958(--26) & 1.731(--25) \\
& 24 & 7.648(--12) & 1.478(--11) & 2.065(--25) & 3.985(--25) \\
& 48 & 2.076(--11) & 4.088(--11) & 5.683(--25) & 1.121(--24) \\
& 98 & 6.042(--11) & 1.201(--10) & 1.504(--24) & 2.987(--24)
\end{tabular}
\end{center}
\medskip

\begin{center}
TABLE 5.1.  Comparison between modified Chebyshev algorithm \\
and modification algorithm in Example 5.1 (cf. Example 3.2)
\end{center}
\bigskip

\noindent
modification, named {\tt test9}, of {\tt test2}. Selected results from 
it, showing the relative discrepancies between
the single-precision values $\alpha_k$, $\beta_k$, resp.
double-precision values
$\alpha_k^d$, $\beta_k^d$, computed by the modified Chebyshev
algorithm and the modification algorithm, are shown in Table
5.1 (cf. Table 3.3).
The maximum errors occur consistently for the last value of
$k$ (= 98).

\bigskip

{\it Example} 5.2.
Induced orthogonal polynomials.
\smallskip

Given an orthogonal polynomial $\pi_m ( \cdot ~; d \lambda )$ of fixed
degree $m \geq 1$, the sequence of orthogonal polynomials
$\hat{\pi}_{k,m} ( \cdot ) = \pi_k ( \cdot ~; \pi_m^2 d \lambda )$,
$k = 0, 1, 2,\ldots,$ has been termed {\it induced orthogonal polynomials}
in [33].
Since their measure $d \hat{\lambda}_m$ modifies the measure
$d \lambda$ by a product of quadratic factors,
$$
d \hat{\lambda}_m (t) = \prod_{\mu =1}^ m
(t - x_{\mu} )^2 \cdot d \lambda (t) ,
\eqno(5.2)
$$
where $x_{\mu}$ are the zeros of $\pi_m$, we can apply the routine
{\tt chri} (with {\tt iopt}=7) $m$ times to compute the $n$
recursion coefficients $\hat{\alpha}_{k,m}$ = $\alpha_k ( d \hat{\lambda}_m )$,
$\hat{\beta}_{k,m}$ = $\beta_k ( d \hat{\lambda}_m )$,
$k = 0, 1,\ldots, n-1$, from the $n+m$
coefficients $\alpha_k = \alpha_k ( d \lambda )$,
$\beta_k = \beta_k ( d \lambda )$, $k = 0, 1,\ldots, n-1+m$.
The subroutines {\tt indp} and {\tt dindp} in the driver {\tt test10}
implement this procedure in single resp.
double precision.
The driver itself uses them to compute the first $n$ = 20 recursion
coefficients of the induced Legendre polynomials with
$m = 0, 1,\ldots, 11$.
It also computes the maximum absolute errors in the
$\hat{\alpha}$'s $( \hat{\alpha}_{k,m} = 0$ for all $m$) and the maximum
relative errors in the $\hat{\beta}$'s by comparing single-precision with
double-precision results.

An excerpt of the output of {\tt test10} is shown in Table 5.2.
It already suggests a high degree of stability of the procedure employed by
{\tt indp}.
This is reinforced by an additional test (not shown in the
package) generating $n$ = 320 recursion
coefficients
$\hat{\alpha}_{k,m}$, $\hat{\beta}_{k,m}$,
$0 \leq k \leq 319$, for $m$ = 40, 80, 160, 320 and
$d \lambda$ being the Legendre, the 1st-kind Chebyshev, the Laguerre,
and the Hermite measure.
Table 5.3 shows the maximum absolute error in the
$\hat{\alpha}_{k,m}$, $0 \leq k \leq 319$ (relative error in the Laguerre
case), and the maximum relative error in the
$\hat{\beta}_{k,m}$, $0 \leq k \leq 319$.

\bigskip
\begin{center}
\begin{tabular}{lrrrr}
\multicolumn{1}{c}{$k$} &
\multicolumn{1}{c}{$m=0$, $\hat{\beta}_{k,m}$} &
\multicolumn{1}{c}{$m=2$, $\hat{\beta}_{k,m}$} &
\multicolumn{1}{c}{$m=6$, $\hat{\beta}_{k,m}$} &
\multicolumn{1}{c}{$m=11$, $\hat{\beta}_{k,m}$} \\ \hline
  & & & &  \\*[-10pt]
\ \, 0 & 2.0000000000 & .1777777778 & .0007380787 & .0000007329 \\
\ \, 1 & .3333333333 & .5238095238 & .5030303030 & .5009523810 \\
\ \, 6 & .2517482517 & .1650550769 & .2947959861 & .2509913424 \\
\, 12 & .2504347826 & .2467060415 & .2521022519 & .1111727541 \\
\, 19 & .2501732502 & .2214990335 & .2274818789 & .2509466619 \\
\\
err $\hat{\alpha}$ & 0.000(0) \ \ \ & 1.350(--13) & 9.450(--13) & 1.357(--12) \\
err $\hat{\beta}$ & 1.737(--14) & 2.032(--13) & 2.055(--12) & 3.748(--12)
\end{tabular}
\end{center}
\medskip

\begin{center}
TABLE 5.2.  Induced Legendre polynomials
\end{center}
\bigskip

\small
\begin{center}
\begin{tabular}{lllllllll}
& \multicolumn{2}{c}{Legendre} &
\multicolumn{2}{c}{Chebyshev} &
\multicolumn{2}{c}{Laguerre} &
\multicolumn{2}{c}{Hermite} \\
$m$ & \ \ \ err $\hat{\alpha}$ & \ \ \ err $\hat{\beta}$ & \ \ \ err $\hat{\alpha}$ & \ \ \ err $\hat{\beta}$ & \ \ \ err $\hat{\alpha}$ & \ \ \ err $\hat{\beta}$ & \ \ \ err $\hat{\alpha}$ & \ \ \ err $\hat{\beta}$ \\ \hline
 & & & & & & & &  \\*[-10pt]
\ \ 40 & 3.4(--11) & 1.5(--10) & 1.9(--9) & 7.9(--10) & 3.0(--10) & 6.0(--10) & 1.8(--9) & 2.7(--10) \\
\ \ 80 & 1.4(--10) & 5.4(--10) & 2.1(--9) & 2.2(--9) & 5.8(--10) & 9.2(--10) & 7.9(--9) & 9.2(--10) \\
160 & 1.5(--9) & 5.1(--9) & 9.5(--9) & 1.1(--8) & 7.8(--10) & 1.4(--9) & 1.1(--8) & 6.8(--10) \\
320 & 3.3(--9) & 2.1(--8) & 9.6(--9) & 2.1(--8) & 3.9(--9) & 7.2(--9) & 2.5(--8) & 1.1(--9)
\end{tabular}
\end{center}
\medskip
\normalsize

\begin{center}
TABLE 5.3.  Accuracy of the recursion coefficients for \\
some classical induced polynomials
\end{center}
\bigskip

5.2.  {\it Methods based on the modified Chebyshev algorithm}.
As was noted earlier, the procedure {\tt chri} becomes unstable for
modified measures involving division of
$d \lambda (t)$ by $t - x$ or
$(t-x)^2 + y^2$ as
$z = x+iy \in \mbox{\myfont C}$ moves away from the ``support interval''
of $d \lambda$, i.e., from the smallest interval containing
the support of $d \lambda$.
We now develop a procedure that works better the further away $z$
is from that interval.

The idea is to use modified moments of $d \hat{\lambda}$ relative to the
polynomials $\pi_k ( \cdot ~; d \lambda )$ to generate the desired
recursion coefficients $\hat{\alpha}_k$,
$\hat{\beta}_k$ via the modified Chebyshev algorithm (cf. \S3).
The modified moments in question are
$$
\nu_k = \nu_k ( x ; d \lambda ) = \int_{\mbox{\myfont R} } ~
\frac{\pi_k ( t ; d \lambda )}{ t - x } ~ d \lambda (t) , ~~~~
k = 0, 1, 2,\ldots,
\eqno(5.3)
$$
for linear divisors, and
$$
\nu_k = \nu_k ( x,y;d \lambda ) = \int_{\mbox{\myfont R} } ~
\frac{\pi_k ( t; d \lambda )}{(t-x)^2 + y^2} ~
d \lambda (t) , ~~~~
k = 0, 1, 2,\ldots,
\eqno(5.4)
$$
for quadratic divisors.
Both can be expressed in terms of the integrals
$$
\begin{array}{r}
{\displaystyle
\rho_k = \rho_k ( z;d \lambda ) = \int_{\mbox{\myfont R}} ~
\frac{\pi_k (t;d \lambda )}{z-t} ~
d \lambda (t) , ~~~~z \in \mbox{\myfont C} \backslash \mbox{supp}
( d \lambda ) , }\
{\displaystyle
k = 0, 1, 2,\ldots, }
\end{array}
\eqno(5.5)
$$
the first by means of
$$
\nu_k ( x; d \lambda ) = - \rho_k
(z ; d \lambda ) , ~~~~ z = x ,
\eqno(5.6)
$$
the others by means of
$$
\nu_k ( x,y;d \lambda ) = - ~
\frac{\mbox{Im} ~\rho_k (z;d \lambda )}{\mbox{Im} ~z} ~, ~~~~
z = x + iy .
\eqno(5.7)
$$

The point to observe is that
$\{ \rho_k ( z ; d \lambda ) \}$ is a minimal solution of the basic
recurrence relation (1.3) for the orthogonal polynomials
$\{ \pi_k ( \cdot ~; d \lambda ) \}$ (cf. [18]).
The quantities $\rho_k ( z ; d \lambda )$,
$k = 0, 1, \ldots , n$, therefore can be computed
accurately by a backward recurrence algorithm ([18, \S5]) which, for
$\nu > n$, produces approximations
$\rho_k^{[ \nu ]} (z;d \lambda )$ converging to
$\rho_k ( z; d \lambda )$ when $\nu \rightarrow \infty$, for any fixed $k$,
$$
\rho_k^{[ \nu ]} ( z ; d \lambda ) \rightarrow
\rho_k ( z ; d \lambda ) , ~~~~ \nu \rightarrow \infty .
\eqno(5.8)
$$
The procedure is implemented in the routine
\medskip

\begin{center}
{\tt knum(n,nu0,numax,z,eps,a,b,rho,nu,ierr,rold)},
\end{center}
\medskip

\noindent
which computes $\rho_k ( z;d \lambda )$ for
$k = 0, 1,\ldots,$ {\tt n} to a relative
precision {\tt eps}.
The results are stored as {\tt rho}$(k) = \rho_{k-1} ( z ; d \lambda )$,
$k = 1, 2,\ldots,$ {\tt n}+1, in the complex array {\tt rho} of dimension
$n+1$.
The user has to provide a starting index {\tt nu0} = $\nu_0 > n$ for
the backward recursion, which the routine then increments by units of 5
until convergence to within {\tt eps} is achieved. If the requested
accuracy {\tt eps} cannot be realized for some $\nu \leq$ {\tt numax},
the routine exits with {\tt ierr} = {\tt numax}. Likewise, if $\nu_0 > $ 
{\tt numax}, the routine exits immediately, with the error flag
{\tt ierr} set equal to {\tt nu0}. 
Otherwise, the value of $\nu$ for which convergence is obtained is
returned as output variable {\tt nu}.
The arrays {\tt a}, {\tt b} of dimension {\tt numax} are to hold the
recursion coefficients
{\tt a}$(k) = \alpha_{k-1} (d \lambda )$,
{\tt b}$(k) = \beta_{k-1} ( d \lambda )$, $k = 1, 2,\ldots,$ {\tt numax,
}
for the given measure $d \lambda$.
The complex array {\tt rold} of dimension $n+1$ is used for
working space.
In the interest of rapid convergence, the routine should be provided
with a realistic estimate of $\nu_0$.
For classical measures, such estimates are known (cf. [18, \S5]) and are
implemented here by the function routines
\medskip

\begin{center}
{\tt nu0jac(n,z,eps), nu0lag(n,z,al,eps), nu0her(n,z,eps)}.
\end{center}
\medskip

\noindent
The first is for Jacobi measures, the second for generalized
Laguerre measures with parameter {\tt al} = $\alpha > -1$, and the last
for the Hermite measure.
Note that $\nu_0$ for Jacobi measures does not depend on the weight
parameters $\alpha$, $\beta$, in contrast to $\nu_0$ for the generalized
Laguerre measure.

The name {\tt knum} comes from the fact that
$\rho_n ( z ; d \lambda )$ in (5.5) is the numerator in the kernel
$$
K_n ( z;d \lambda ) = ~
\frac{ \rho_n ( z;d \lambda )}{\pi_n (z;d \lambda )}
\eqno(5.9)
$$
of the remainder term of the $n$-point Gaussian quadrature
rule for analytic functions (cf., e.g., [36]).
For the sequence of kernels $K_0$, $K_1,\ldots, K_n$ we have the following
routine:
\medskip

\begin{verbatim}
     subroutine kern(n,nu0,numax,z,eps,a,b,ker,nu,ierr,rold)
     complex z,ker,rold,p0,p,pm1
     dimension a(numax),b(numax),ker(*),rold(*)
     call knum(n,nu0,numax,z,eps,a,b,ker,nu,ierr,rold)
     if(ierr.ne.0) return
     p0=(0.,0.)
     p=(1.,0.)
     do 10 k=1,n
       pm1=p0
       p0=p
       p=(z-a(k))*p0-b(k)*pm1
       ker(k+1)=ker(k+1)/p
  10 continue
     return
     end
\end{verbatim}
\medskip

\noindent
The meaning of the input and output parameters is the same as in {\tt knum}.
The double precision version of the routine is named {\tt
dkern}.

All the ingredients are now in place to describe the workings of
{\tt gchri}, the alternative routine to {\tt chri} when the latter is
unstable.
First, the routine {\tt knum} is used to produce the first $2n$ modified
moments $\nu_k ( x;d \lambda )$ resp. $\nu_k (x,y;d \lambda )$,
$k = 0, 1,\ldots, 2n-1$.  These are then supplied
to the routine {\tt cheb} along with the recursion coefficients
$\alpha_k ( d \lambda )$, $\beta_k ( d \lambda )$
(needed anyhow for the computation of the $\nu_k$), which produces
the desired coefficients $\alpha_k ( d \hat{\lambda} )$,
$\beta_k ( d \hat{\lambda} )$, $k = 0, 1,\ldots, n-1$.
The routine has the calling sequence
\medskip

\noindent
\begin{verbatim}
     gchri(n,iopt,nu0,numax,eps,a,b,x,y,alpha,beta,
        nu,ierr,ierrc,fnu,rho,rold,s,s0,s1,s2).
\end{verbatim}
\medskip

\noindent
\begin{tabular}{lp{3in}}
On entry,  & \\
\\
\ \ \ \ \ {\tt n} & is the number of recursion coefficients
desired;
type integer \\
\\
\ \ \ \ \ {\tt iopt} & is an integer identifying the type of modification
as follows: \\
& 1: \ \ $d \hat{\lambda} (t) = d \lambda (t) / ( t-x )$, where $x$ is assumed \\
& \hspace{.25in}  outside the smallest interval containing \\
& \hspace{.25in}  supp$(d \lambda )$ \\
& 2: \ \ $d \hat{\lambda} (t) = d \lambda (t) / ((t-x)^2 + y^2 )$, $y > 0$ \\
\\
\ \ \ \ \ {\tt nu0} & is an integer $\nu_0 \geq 2n$ estimating the starting
index for the backward recursion to compute the modified moments; if no
other choices are available, take {\tt nu0} = $3\times${\tt n} 
\end{tabular}

\noindent
\begin{tabular}{lp{3in}}
\\
\ \ \ \ \ {\tt numax} & is an integer used to terminate backward
recursion in case of nonconvergence; a conservative choice is
{\tt numax}=500 \\
\\
\ \ \ \ \ {\tt eps} & is a relative error tolerance; type real \\
\\
\ \ \ \ \ {\tt a,b} & are arrays of dimension {\tt numax} to be supplied
with the recursion coefficients {\tt a}$(k) = \alpha_{k-1} ( d \lambda )$,
{\tt b}$(k) = \beta_{k-1} ( d \lambda )$, $k = 1, 2,\ldots,$
{\tt numax} \\
\\
\ \ \ \ \ {\tt x,y} & are real parameters defining the linear
and quadratic divisors of $d \lambda$. \\
\\
On return, & \\
\\
\ \ \ \ \ {\tt alpha,beta} & are arrays of dimension {\tt n} containing
the desired recursion coefficients {\tt alpha}$(k) = \hat{\alpha}_{k-1}$
{\tt beta}$(k) = \hat{\beta}_{k-1}$, $k = 1, 2,\ldots,$ {\tt n} \\
\\
\ \ \ \ \ {\tt nu} & is the index $\nu$ for which the error tolerance
{\tt eps} is satisfied for the first time; if it is never satisfied,
{\tt nu} will have the value {\tt numax} \\
\\
\ \ \ \ \ {ierr} & is an error flag, where \\
& {\tt ierr} = 0 on normal return \\
& {\tt ierr} = 1 if {\tt iopt} is inadmissible \\
& {\tt ierr} = {\tt nu0} if {\tt nu0} $>$ {\tt numax} \\
& {\tt ierr} = {\tt numax} if the backward recurrence \\
& \hspace{.5in} algorithm does not converge. \\
& {\tt ierr} = --1 if {\tt n} is not in range
\\
\ \ \ \ \ {\tt ierrc} & is an error flag inherited from the
routine {\tt cheb}.
\end{tabular}
\medskip

\noindent
The real arrays {\tt fnu,s,s0,s1,s2} are working space,
all of dimension 2$\times${\tt n} except {\tt s}, which has dimension
{\tt n}.
The complex arrays {\tt rho}, {\tt rold} are also working space, both
of dimension $2n$.
The routine calls on the subroutines {\tt knum} and {\tt cheb}.
The double-precision version of {\tt gchri} has the name
{\tt dgchri}.

Since the routine {\tt gchri} is based on the modified Chebyshev
algorithm, it shares with the latter its proneness
to ill-conditioning, particularly in cases
of measures supported on an infinite interval.
On finitely supported measures, however, it can be quite effective,
as will be seen in the next example.
\bigskip

{\it Example} 5.3.
The performance of {\tt chri} and {\tt gchri}.
\smallskip

To illustrate the severe limitations of the routine
{\tt chri} in situations where divisions of the measure $d \lambda$ are
involved, and at the same time to document the effectiveness
of {\tt gchri}, we ran both routines with $n$ = 40 for Jacobi
measures
$d \lambda^{( \alpha , \beta )}$ with parameters
$\alpha , \beta = - .8(.4).8$,
$\beta \geq \alpha$.
This is done in {\tt test11}.

The routine {\tt test11} first tests division by
$t - x$, where
$x$ = --1.001,\linebreak  --1.01 , --1.04, --1.07, --1.1.
Both routines {\tt chri} and {\tt gchri} are run in single and double
precision, the latter with $\varepsilon = 10\times \varepsilon^{s}$ and
$\varepsilon = 100\times \varepsilon^{d}$, respectively.
The double-precision results are used to determine the absolute
errors in the $\hat{\alpha}$'s and the relative errors in the
$\hat{\beta}$'s for each routine.
The required coefficients
$\alpha_k$, $\beta_k$, $0 \leq k \leq \nu_{\max} - 1$
$( \nu_{\max} = 500$ for single precision, and 800 for double precision)
are supplied by {\tt recur} and {\tt drecur} with
{\tt ipoly} = 6.
The routine {\tt nu0jac} is used to provide the starting recurrence
index $\nu_0$ resp. $\nu_0^d$.
In Tables 5.4 and 5.5, relating respectively to linear and quadratic
divisors, we give only the results for the Legendre measure
($ \alpha = \beta = 0)$.
The first line in each 3-line block of Table 5.4 shows $x$,
$\nu_0$, $\nu_0^d$ and the maximum (over $k$,
$0 \leq k \leq 39$) errors in the
$\hat{\alpha}_k$ and $\hat{\beta}_k$ for {\tt gchri}, followed by the
analogous information (excepting the $\nu_0$'s) for {\tt chri}.
The recurrence index $\nu$ yielding convergence was found (not
shown in {\tt test11}) to be
$\nu = \nu_0 + 5$ and
$\nu^d = \nu_0^d + 5$, without exception.
\bigskip

\noindent
\begin{tabular}{lrcrrll}
& & &
\multicolumn{2}{c}{\tt gchri} &
\multicolumn{2}{c}{\tt chri} \\
\multicolumn{1}{c}{$x$} &
\multicolumn{1}{c}{$\nu_0$} &
\multicolumn{1}{c}{$\nu_{0}^{d}$} &
\multicolumn{1}{c}{err $\hat{\alpha}$} &
\multicolumn{1}{c}{err $\hat{\beta}$} &
\multicolumn{1}{c}{err $\hat{\alpha}$} &
\multicolumn{1}{c}{err $\hat{\beta}$} \\ \hline
 & & & & &  \\*[-10pt]
--1.001 & 418 & 757 & 8.000(--14) & 1.559(--13) & 1.013(--13) & 1.647(--13) \\
& \multicolumn{2}{c}{reconstruction} & 8.527(--14) &
1.705(--13) & 1.010(--27) & 2.423(--27) \\
& \multicolumn{2}{c}{errors} & 1.421(--14) & 5.329(--14) &
2.019(--28) & 1.211(--27) \\
\\
--1.01 & 187 & 294 & 4.016(--14) & 6.907(--14) & 1.396(--10) &
2.424(--10) \\
& & & 3.553(--14) & 9.946(--14) & 6.058(--28) & 1.211(--27) \\
& & & 7.105(--15) & 4.262(--14) & 1.515(--28) & 9.080(--28) \\
\\
--1.04 & 133 & 187 & 3.590(--14) & 4.759(--14) & 5.944(--6) &
8.970(--6) \\
& & & 2.842(--14) & 7.103(--14) & 5.554(--28) & 1.312(--27) \\
& & & 7.105(--15) & 4.263(--14) & 1.010(--28) & 9.080(--28) \\
\\
--1.07 & 120 & 161 & 2.194(--14) & 4.850(--14) & 5.334(--3) &
7.460(--3) \\
& & & 2.842(--14) & 7.104(--14) & 6.058(--28) & 1.211(--27) \\
& & & 7.105(--15) & 4.263(--14) & 1.010(--28) & 7.062(--28) \\
\\ 
--1.1 & 114 & 148 & 2.238(--14) & 4.359(--14) & 4.163(0) &
4.959(+1) \\
& & & 2.132(--14) & 5.683(--14) & 3.534(--28) & 1.009(--27) \\
& & & 1.549(--12) & 1.833(--12) & 1.010(--28) & 6.057(--28) \\
\end{tabular}
\medskip

\begin{center}
TABLE 5.4.  Performance of {\tt gchri} and {\tt chri} for elementary \\
divisors $t-x$ of the Legendre measure $d \lambda (t)$
\end{center}
\bigskip

It can be seen from the leading lines in Table 5.4 that {\tt chri}
rapidly loses accuracy as $x$ moves away from the interval [$-$1,1],
all single-precision accuracy being gone by the time $x$
reaches --1.1.
Similar, if not more rapid, erosion of accuracy is observed for the other
parameter values of $\alpha , \beta$.
The next two lines in each 3-line block show ``reconstruction errors'',
i.e., the maximum errors in the $\alpha$'s and $\beta$'s if the
$\hat{\alpha}$'s and $\hat{\beta}$'s produced by
{\tt gchri}, {\tt chri} and {\tt dgchri}, {\tt dchri}
are fed back to the routines {\tt chri} and {\tt dchri} with
{\tt iopt} = 1 to recover the original recursion coefficients in single and
double precision.
The first of these two lines shows the errors in reconstructing these
coefficients from the output of {\tt gchri} (resp. {\tt dgchri}),
the second from the output of {\tt chri} (resp. {\tt dchri}).
Rather remarkably, the coefficients are recovered to essentially full accuracy,
even when the input coefficients (produced by {\tt chri} and
{\tt dchri}) are very inaccurate!
This is certainly a phenomenon that deserves further study.
It can also be seen from Table 5.4 (and the more complete results in \S1
of the package) that {\tt gchri} consistently produces accurate results,
some slight deterioration occurring only very close to $x = -1$, where the
routine has to work harder.

The second half of {\tt test11}
tests division by
$(t-x)^2 + y^2$, where
$z = x+iy$ is taken along the upper half of the ellipse
$$
{\cal E}_{\rho} = \{ z \in \mbox{\myfont C}: ~
z = \frac{1}{2} ~ \left( \rho e^{i \vartheta} +
\frac{1}{\rho} e^{- i \vartheta} \right) , ~
0 \leq \vartheta \leq 2 \pi \}, ~ \rho > 1,
\eqno(5.10)
$$
which has foci $\pm 1$ and sum of the semiaxes equal to $\rho$.
(These ellipses are contours of constant $\nu_0$ for Jacobi measures.)
We generated information analogous to the one in Table 5.4, for
$\rho$ = 1.05, 1.1625, 1.275, 1.3875, 1.5, except that all quantities
are averaged over 19 equally spaced points on
${\cal E}_{\rho}$ corresponding to
$\vartheta = j \pi / 20$,
$j = 1, 2,\ldots, 19$.
Selected results (bars indicate averaging),
again for the Legendre case, are shown in Table 5.5.
They reveal a behavior very similar to the one in Table 5.4 for linear
divisors.

\bigskip

\noindent
\begin{tabular}{lrcllll}
& & &
\multicolumn{2}{c}{\tt gchri} &
\multicolumn{2}{c}{\tt chri} \\
\multicolumn{1}{c}{$\rho$} &
\multicolumn{1}{c}{$\bar{\nu}_0$} &
\multicolumn{1}{c}{$\bar{\nu}_{0}^{d}$} &
\multicolumn{1}{c}{\( \overline{err} \) $\hat{\alpha}$} &
\multicolumn{1}{c}{\( \overline{err} \) $\hat{\beta}$} &
\multicolumn{1}{c}{\( \overline{err} \) $\hat{\alpha}$} &
\multicolumn{1}{c}{\( \overline{err} \) $\hat{\beta}$} \\ \hline
 & & & & &  \\*[-10pt]
1.05 & 390 & 700 & 7.879(--13) & 1.440(--12) & 7.685(--14) & 1.556(--13) \\
& \multicolumn{2}{c}{reconstruction} & 7.814(--13) &
1.433(--12) & 1.768(--26) & 3.042(--26) \\
& \multicolumn{2}{c}{errors} & 2.024(--14) & 8.442(--14) &
3.016(--28) & 1.742(--27) \\
\\
1.275 & 142 & 204 & 6.252(--14) & 1.287(--13) & 4.562(--7) &
6.162(--7) \\
& & & 6.554(--14) & 1.279(--13) & 1.541(--27) & 3.061(--27) \\
& & & 2.295(--14) & 8.970(--14) & 3.579(--28) & 1.646(--27) \\
\\
1.5 & 117 & 154 & 3.991(--14) & 7.966(--14) & 4.906(--1) &
2.339(0) \\
& & & 4.207(--14) & 9.064(--14) & 6.932(--28) & 1.676(--27) \\
& & & 3.805(--14) & 8.971(--14) & 4.351(--28) & 1.744(--27)
\end{tabular}
\medskip

\begin{center}
TABLE 5.5.  Performance of {\tt gchri} and {\tt chri} for elementary \\
divisors $(t-x)^2 + y^2$ of the Legendre measure $d \lambda (t)$ with $z = x+iy$ on ${\cal E}_{\rho}$.
\end{center}

\newpage

\noindent
{\small 6.  GAUSS-TYPE QUADRATURE RULES}

\bigskip

One of the important uses of orthogonal polynomials is in the
approximation of integrals involving a positive measure
$d \lambda$ by quadrature rules of maximum, or nearly maximum,
algebraic degree of exactness.
In this context, it is indispensable to know the recursion coefficients
for the respective orthogonal polynomials
$\{ \pi_k ( \cdot ~; d \lambda ) \}$, since they allow us to generate
the desired quadrature rules accurately and effectively via eigenvalue
techniques.
The software developed in the previous sections thus finds here a vast
area of application.
\bigskip

6.1. {\it Gaussian quadrature}.
Given the (positive) measure $d \lambda$ (having an infinite number of
support points),
there exists, for each $n \in \mbox{\myfont N}$, a quadrature rule
$$
\int_{\mbox{\myfont R}} f (t) d \lambda (t) =
\sum_{k=1}^n w_k f (x_k ) + R_n (f)
\eqno(6.1)
$$
having algebraic degree of exactness
$2n - 1$, i.e., zero error, $R_n (f) = 0$, whenever $f$ is a polynomial of
degree $\leq 2n-1$.
The nodes $x_k$ indeed are the zeros of the $n$th-degree
orthogonal polynomial
$\pi_n ( \cdot ~; d \lambda )$, and the weights $w_k$, which are all
positive, are also expressible in terms of the same orthogonal polynomials.
Alternatively, and more importantly for computational purposes,
the nodes $x_k$ are the eigenvalues of the $n$th-order Jacobi matrix
$$
J_n ( d \lambda ) = \left[ \begin{array}{ccccc}
\alpha_0 & \sqrt{\beta_1} & & & 0 \\
\sqrt{\beta_1} & \alpha_1 & \sqrt{\beta_2} & & \\
& \sqrt{\beta_2} & \ddots & \ddots & \\
& & \ddots & & \sqrt{\beta_{n-1} } \\
0 & & & \sqrt{\beta_{n-1}} & \alpha_{n-1} \\
\end{array}
\right] ,
\eqno(6.2)
$$
where $\alpha_k = \alpha_k ( d \lambda )$,
$\beta_k = \beta_k ( d \lambda )$ are the recurrence
coefficients for the (monic) orthogonal polynomials
$\{ \pi_k ( \cdot ~; d \lambda ) \}$, and the weights
$w_k$ are expressible in terms of the associated eigenvectors.
Specifically, if
$$
J_n ( d \lambda ) v_k = x_k v_k , ~~~~
v_k^T v_k = 1, ~~~~
k = 1, 2,\ldots, n,
\eqno(6.3)
$$
i.e., $v_k$ is the normalized eigenvector of
$J_n (d \lambda )$ corresponding to the eigenvalue $x_k$, then
$$
w_k = \beta_0 v_{k,1}^2 , ~~~~
k = 1, 2,\ldots,n,
\eqno(6.4)
$$
where $\beta_0 = \beta_0 ( d \lambda )$ is defined in (1.4) and
$v_{k,1}$ is the first component of $v_k$
(cf. [40]).
There are well-known and efficient algorithms, such as the $QR$ algorithm,
to compute eigenvalues and (part of the) eigenvectors of symmetric
tridiagonal matrices.
These are used in the routine {\tt gauss}\footnote{This routine was
kindly supplied to the author by Professor G.H. Golub.}, whose
calling sequence is
\medskip

\begin{center}
{\tt gauss(n,alpha,beta,eps,zero,weight,ierr,e).}
\end{center}
\bigskip

\noindent
\begin{tabular}{lp{3in}}
On entry,  & \\
\\
\ \ \ \ \ {\tt n} & is the number of terms in the Gauss formula;
type integer \\
\\
\ \ \ \ \ {\tt alpha,beta} & are arrays of dimension {\tt n} assumed
to hold the recursion coefficients {\tt alpha}$(k) = \alpha_{k-1}$,
{\tt beta}$(k) = \beta_{k-1}$, $k = 1, 2,\ldots,$ {\tt n}
\\
\ \ \ \ \ {\tt eps} & is a relative error tolerance, for example, the machine
precision. \\
\\
On return, & \\
\\
\ \ \ \ \ {\tt zero,weight} & are arrays of dimension {\tt n} containing
the nodes (in increasing order) and the corresponding weights of the Gauss
formula, {\tt zero}$(k) = x_k$, {\tt weight}$(k) = w_k$,
$k = 1, 2,\ldots,$ {\tt n} \\
\\
\ \ \ \ \ {\tt ierr} & is an error flag equal to 0 on normal return, 
equal to $i$ if the $QR$ algorithm does not converge 
within 30 iterations on evaluating the $i$th eigenvalue, equal to --1
if {\tt n} is not in range, and equal to --2 if one of the $\beta$'s
is negative.
\end{tabular}
\medskip

\noindent
The array {\tt e} of dimension {\tt n} is used for working space.
The double precision routine has the name {\tt dgauss}.

We refrain here from giving numerical examples,
since the use of the routine {\tt gauss},
and the routines yet to be described, is straightforward.
Some use of {\tt gauss} and {\tt dgauss} has already been made in
Examples 4.2, 4.3, 4.4 and 5.2.

\bigskip
6.2. {\it Gauss-Radau quadrature}.
We now assume that $d \lambda$ is a measure whose support is either bounded
from below, or bounded from above, or both.
Let, then, $x_0$ be either the infimum or the supremum of
supp $d \lambda$, so long as it is finite.
(Typically, if supp $d \lambda = [-1,1]$, then $x_0$ could be either --1
or +1; if supp $d \lambda = [0, \infty ]$, then $x_0$ would have to be 0;
etc.).
By {\it Gauss-Radau quadrature} we then mean a quadrature rule of maximum
degree of exactness that contains among the nodes the point $x_0$.
It thus has the form
$$
\int_{\mbox{\myfont R}} f(t) d \lambda (t) =
w_0 f(x_0) + \sum_{k=1}^n w_k f (x_k ) + R_n (f) ,
\eqno(6.5)
$$
and, as is well known, can be made to have degree of exactness
$2n$, i.e., $R_n (f) = 0$
for all polynomials of degree $\leq 2n$.
Interestingly, all nodes
$x_0$, $x_1,\ldots, x_n$ and weights
$w_0$, $w_1 ,\ldots, w_n$ can again be interpreted in terms of
eigenvalues and eigenvectors, exactly as in the case of
Gaussian quadrature rules, but now relative to the matrix (cf. [39])

$$
J_{n+1}^* (d \lambda ) = \left[ \begin{array}{ccccc}
\alpha_0 & \sqrt{\beta_1} & & & 0 \\
\sqrt{\beta_1} & \alpha_1 & \ddots &  \\
& \ddots & \ddots & \sqrt{\beta_{n-1}} & \\
& & \sqrt{\beta_{n-1} } & \alpha_{n-1} & \sqrt{\beta_n} \\
0 & & & \sqrt{\beta_{n} } & \alpha_n^*
\end{array}
\right] ~ \in ~ \mbox{\myfont R}^{(n+1) \times (n+1) } ,
\eqno(6.6)
$$
where $\alpha_k = \alpha_k ( d \lambda )$
$( 0 \leq k \leq n-1)$,
$\beta_k = \beta_k ( d \lambda )$
$( 1 \leq k \leq n )$ as
before, but
$$
\alpha_n^* = \alpha_n^* ( d \lambda ) = x_0 - \beta_n ~
\frac{ \pi_{n-1} ( x_0 ; d \lambda )}{\pi_n ( x_0 ; d \lambda )} ~.
\eqno(6.7)
$$
Hence, we can use the routine {\tt gauss}
to generate the Gauss-Radau formula.
This is done in the following subroutine.

\medskip
\begin{verbatim}
      subroutine radau(n,alpha,beta,end,zero,weight,ierr,e,a,b)
      dimension alpha(*),beta(*),zero(*),weight(*),e(*),a(*),b(*)
c
c The arrays  alpha,beta,zero,weight,e,a,b  are assumed to have
c dimension n+1.
c
      epsma=r1mach(3)
c
c epsma is the machine single precision.
c
      np1=n+1
      do 10 k=1,np1
        a(k)=alpha(k)
        b(k)=beta(k)
   10 continue
      p0=0.
      p1=1.
      do 20 k=1,n
        pm1=p0
        p0=p1
        p1=(end-a(k))*p0-b(k)*pm1
   20 continue
      a(np1)=end-b(np1)*p0/p1
      call gauss(np1,a,b,epsma,zero,weight,ierr,e)
      return
      end
\end{verbatim}
\medskip

\noindent
The input variables are {\tt n}, {\tt alpha},
{\tt beta}, {\tt end} representing, respectively, $n$,
two arrays of dimension $n + 1$ containing the
$\alpha_k ( d \lambda )$, $\beta_k ( d \lambda )$,
$k = 0, 1, 2,\ldots, n$, and the ``endpoint''
$x_0$.
The nodes (in increasing order) of the Gauss-Radau formula
are returned in the array {\tt zero}, the corresponding weights
in the array {\tt weight}.
The arrays {\tt e}, {\tt a}, {\tt b} are working space, and
{\tt ierr} an error flag inherited from the routine {\tt gauss}.
The double-precision routine has the name {\tt dradau}.

We remark that $x_0$ could also be outside the support of
$d \lambda$, in which case the routine would generate a
``Christoffel-type'' quadrature rule.
\bigskip

6.3.  {\tt Gauss-Lobatto quadrature}.
Assuming now the support of $d \lambda$ bounded on either side,
we let $x_0$ = inf supp $( d \lambda )$ and
$x_{n+1}$ = sup supp $( d \lambda )$ and consider a quadrature
rule of the type
$$
\int_{\mbox{\myfont R}} f(t) d \lambda (t) =
w_0 f(x_0 ) + \sum_{k=1}^n w_k f (x_k ) +
w_{n+1} f (x_{n+1}) + R_n (f)
\eqno(6.8)
$$
having maximum degree of exactness $2n+1$.
This is called the {\it Gauss-Lobatto quadrature rule}.
Its nodes $x_0$, $x_1 ,\ldots, x_{n+1}$ and weights
$w_0$, $w_1 ,\ldots, w_{n+1}$ again admit the same spectral
representation as in the case of Gauss and Gauss-Radau formulae, only
this time the matrix in question has order $n+2$ and is given by
(cf. [39])

$$
J_{n+2}^* ( d \lambda ) = \left[ \begin{array}{cccccc}
\alpha_0 & \sqrt{\beta_1} & & & & 0 \\
\sqrt{\beta_1} & \alpha_1 & \sqrt{\beta_2} & & & \\
& \sqrt{\beta_2} & & \ddots & & \\
& & \ddots & \ddots & \sqrt{\beta_n} & \\
& & & \sqrt{\beta_n} & \alpha_n & \sqrt{\beta_{n+1}^*} \\
0 & & & & \sqrt{\beta_{n+1}^*} & \alpha_{n+1}^*    
\end{array}
\right] ~.
\eqno(6.9)
$$
Here, as before, $\alpha_k = \alpha_k ( d \lambda )$
$(0 \leq k \leq n)$, $\beta_k = \beta_k ( d \lambda )$
$(1 \leq k \leq n)$, and 
$\alpha_{n+1}^*$, $\beta_{n+1}^*$ are the unique solution of the
linear $2 \times 2$ system
$$
\left[ \begin{array}{ll} 
\pi_{n+1} (x_0 ; d \lambda ) & \pi_n ( x_0 ; d \lambda ) \\
\pi_{n+1} (x_{n+1} ; d \lambda ) & \pi_n ( x_{n+1} ; d \lambda ) 
\end{array} \right] ~~
\left[ \begin{array}{l}
\alpha_{n+1}^* \\ \beta_{n+1}^* \end{array}
\right] = \left[ \begin{array}{l}
x_0 \pi_{n+1} (x_0 ; d \lambda ) \\
x_{n+1} \pi_{n+1} ( x_{n+1} ; d \lambda ) 
\end{array} \right] ~.
\eqno(6.10)
$$
Hence, we have the following routine for generating Gauss-Lobatto
formulae:
\medskip

\begin{verbatim}
      subroutine lob(n,alpha,beta,aleft,right,zero,weight,ierr,e,a,b)
      dimension alpha(*),beta(*),zero(*),weight(*),e(*),a(*),b(*)
c
c The arrays  alpha,beta,zero,weight,e,a,b  are assumed to have
c dimension n+2.
c
      epsma=r1mach(3)
c
c epsma is the machine single precision.
c
      np1=n+1 
      np2=n+2
      do 10 k=1,np2
        a(k)=alpha(k)
        b(k)=beta(k)
   10 continue
      p0l=0.
      p0r=0.
      p1l=1.
      p1r=1.
      do 20 k=1,np1
        pm1l=p0l
        p0l=p1l
        pm1r=p0r
        p0r=p1r
        p1l=(aleft-a(k))*p0l-b(k)*pm1l
        p1r=(right-a(k))*p0r-b(k)*pm1r
   20 continue
      det=p1l*p0r-p1r*p0l
      a(np2)=(aleft*p1l*p0r-right*p1r*p0l)/det
      b(np2)=(right-aleft)*p1l*p1r/det
      call gauss(np2,a,b,epsma,zero,weight,ierr,e)
      return
      end
\end{verbatim}
\medskip

\noindent
The meaning of the input and output variables is as in the
routine {\tt radau}, the variable {\tt aleft} standing for
$x_0$ and {\tt right} for $x_{n+1}$. The double-precision
routine is named {\tt dlob}.

A remark analogous to the one after the routine {\tt radau}
applies to the routine {\tt lob}.

\vspace{.5in}

{\tt Acknowledgment.} The author gratefully acknowledges a number
of suggestions from two anonymous referees and from the associate
editor, Dr. Ronald F. Boisvert, for improving the code of the
package.

\vspace{.5in}

\noindent
Walter Gautschi

\noindent
Department of Computer Sciences

\noindent
Purdue University

\noindent
West Lafayette, IN 47907

\noindent
{\it E-mail address}:  wxg@cs.purdue.edu

\newpage

\noindent
{\small REFERENCES}
\medskip
\begin{enumerate}
\item
{\tt ABRAMOWITZ, M.} and {\tt STEGUN, I.A.} (eds.),
Handbook of Mathematical Functions,
{\it NBS Appl. Math. Ser. 55}, U.S. Government Printing
Office, Washington, D.C., 1964.
\item
{\tt BOLEY, D.} and {\tt GOLUB, G.H.},
A survey of matrix inverse eigenvalue problems, 
{\it Inverse Problems 3} (1987), 595--622.
\item
{\tt de BOOR, C.} and {\tt GOLUB, G.H.},
The numerically stable reconstruction of a Jacobi matrix from
spectral data, {\it Linear Algebra Appl. 21} (1978), 245--260.
\item
{\tt CHEBYSHEV, P.L.},
Sur l'interpolation par la m\'{e}thode des moindres
carr\'{e}s, {\it M\'{e}m. Acad. Imp\'{e}r. Sci. St.
Petersbourg (7) 1}, no. 15 (1859), 1--24.
[{\it {\OE}uvres I}, pp. 473--498.]
\item
{\tt CHIHARA, T.S.},
{\it An Introduction to Orthogonal Polynomials},
Gordon and Breach, New York, 1978.
\item
{\tt CHIHARA, T.S.},
Orthogonal polynomials and measures with end point masses,
{\it Rocky Mountain J. Math. 15} (1985), 705--719.
\item
{\tt CHRISTOFFEL, E.B.},
\"{U}ber die Gau$\beta$ische Quadratur und eine Verallgemeinerung
derselben, {\it J. Reine Angew. Math. 55} (1858),
61--82. 
[{\it Ges. Math. Abhandlungen I}, pp. 65--87.]
\item
{\tt CHRISTOFFEL, E.B.},
Sur une classe particuli\`{e}re de fonctions enti\`{e}res
et de fractions continues, {\it Ann. Mat. Pura Appl. (2) 8}
(1877), 1--10.
[{\it Ges. Math. Abhandlungen II}, pp. 42--50.]
\item
{\tt CODY, W.J.} and {\tt HILLSTROM, K.E.},
Chebyshev approximations for the natural logarithm of the gamma
function, {\it Math. Comp. 21} (1967), 198--203.
\item
{\tt DANLOY, B.},
Numerical construction of Gaussian quadrature formulas for
$\int_0^1 (- \mbox{Log} ~x) \cdot x^{\alpha} \cdot f(x) \cdot dx$
and $\int_0^{\infty} E_m (x) \cdot f(x) \cdot dx$,
{\it Math. Comp. 27} (1973), 861--869.
\item  
{\tt FRONTINI, M., GAUTSCHI, W.}                   
and {\tt MILOVANOVI\'{C}, G.V.},
Moment-\linebreak preserving
spline approximation on finite intervals,
{\it Numer. Math. 50} (1987), 503--518.
\item
{\tt GALANT, D.}, Gauss quadrature rules for the
evaluation of\linebreak $2 \pi^{- 1/2} \int_0^{\infty} \exp (-x^2 ) f(x) dx$,
Review 42, {\it Math. Comp. 23} (1969), 676--677.
Loose microfiche suppl. E.
\item
{\tt GALANT, D.},
An implementation of Christoffel's theorem in the theory of
orthogonal polynomials, {\it Math. Comp. 25} (1971),
111--113.
\item
{\tt GALANT, D.},
Algebraic methods for modified orthogonal polynomials, {\it Math.
Comp. 59} (1992), 541--546.
\item
{\tt GAUTSCHI, W.},
Numerical quadrature in the presence of a singularity,
{\it SIAM J. Numer. Anal. 4} (1967), 357--362.
\item
{\tt GAUTSCHI, W.},
Computational aspects of three-term recurrence relations,
{\it SIAM Rev. 9} (1967), 24--82.
\item
{\tt GAUTSCHI, W.},
On the preceding paper ``A Legendre polynomial integral''
by James L. Blue, {\it Math. Comp. 33} (1979), 742--743.
\item
{\tt GAUTSCHI, W.},
Minimal solutions of three-term recurrence relations and
orthogonal polynomials, {\it Math. Comp. 36} (1981),
547--554.
\item
{\tt GAUTSCHI, W.},
On generating orthogonal polynomials,
{\it SIAM J. Sci. Stat. Comput. 3} (1982), 289--317.
\item
{\tt GAUTSCHI, W.},
An algorithmic implementation of the generalized\linebreak Christoffel
theorem, in {\it Numerical Integration} 
(G. H\"{a}mmerlin, ed.), pp. 89--106.
{\it Intern. Ser. Numer. Math.}, v. {\it 57},
Birkh\"{a}user, Basel, 1982.
\item
{\tt GAUTSCHI, W.},
Discrete approximations to spherically symmetric distributions,
{\it Numer. Math. 44} (1984), 473--483.
\item
{\tt GAUTSCHI, W.},
Questions of numerical condition related to polynomials, in
{\it Studies in Mathematics}, v. {\it 24}, {\it Studies in Numerical
Analysis} (G.H. Golub, ed.), pp. 140--177.  The Mathematical Association of
America, 1984.
\item
{\tt GAUTSCHI, W.},
On the sensitivity of orthogonal polynomials to perturbations
in the moments, {\it Numer. Math. 48} (1986), 369--382.
\item
{\tt GAUTSCHI, W.},
Reminiscences of my involvement in de Branges's proof of the
Bieberbach conjecture, in {\it The Bieberbach
Conjecture} (A. Baernstein II et al., eds.), {\it Mathematical
Surveys and Monographs}, No. {\it 21}, pp. 205--211.   American Math. Soc.,
Providence, R.I., 1986.
\item
{\tt GAUTSCHI, W.},
Computational aspects of orthogonal polynomials, in 
{\it Orthogonal Polynomials -- Theory and Practice}
(P. Nevai, ed.), pp. 181--216.
{\it NATO ASI Series, Series C:
Mathematical and Physical Sciences}, v. {\it 294}, Kluwer, 
Dordrecht, 1990.
\item
{\tt GAUTSCHI, W.},
A class of slowly convergent series and their summation by
Gaussian quadrature, {\it Math. Comp. 57} (1991), 309--324.
\item
{\tt GAUTSCHI, W.},
On certain slowly convergent series occurring in plate
contact problems, {\it Math. Comp. 57} (1991), 325--338.
\item
{\tt GAUTSCHI, W.},
Quadrature formulae on half-infinite intervals,
{\it BIT 31} (1991), 438--446.
\item
{\tt GAUTSCHI, W.},
Computational problems and applications of orthogonal 
polynomials, in {\it Orthogonal Polynomials and Their
Applications} (C. Brezinski et al., eds.), pp. 61--71.
{\it IMACS Annals on Computing and Applied Mathematics}, v. {\it 9},
Baltzer, Basel, 1991.
\item
{\tt GAUTSCHI, W.},
Gauss-type quadrature rules for rational functions, in {\it Numerical
Integration IV} (H. Brass and G. H\"{a}mmerlin, eds.), Internat.
Ser. Numer. Math., to appear.
\item
{\tt GAUTSCHI, W.},
On the computation of generalized Fermi-Dirac and Bose-Einstein
integrals,
{\it Comput. Phys. Comm.}, to appear.
\item
{\tt GAUTSCHI, W.},
Is the recurrence relation for orthogonal polynomials always
stable?, {\it BIT}, to appear.
\item
{\tt GAUTSCHI, W.}, and {\tt LI, S.},
A set of orthogonal polynomials induced by a given orthogonal
polynomial, {\it Aequationes Math.}, to appear.
\item
{\tt GAUTSCHI, W.}, and {\tt MILOVANOVI\'{C}, G.V.},
Gaussian quadrature involving Einstein and Fermi functions with an
application to summation of series,
{\it Math. Comp. 44} (1985), 177--190.
\item
{\tt GAUTSCHI, W.},
and {\tt MILOVANOVI\'{C}, G.V.},
Spline approximations to spherically symmetric distributions,
{\it Numer. Math. 49} (1986), 111--121.
\item
{\tt GAUTSCHI, W.} and {\tt VARGA, R.S.},
Error bounds for Gaussian quadrature of analytic functions,
{\it SIAM J. Numer. Anal. 20} (1983), 1170--1186.
\item
{\tt GAUTSCHI, W.} and {\tt WIMP, J.},
Computing the Hilbert transform of a Jacobi weight function,
{\it BIT 27} (1987), 203--215.
\item
{\tt GAUTSCHI, W., KOVA\v{C}EVI\'{C}, M.A.}, and
{\tt MILOVANOVI\'{C}, G.V.},
The numerical evaluation of singular integrals with 
coth-kernel, {\it BIT 27} (1987), 389--402.
\item
{\tt GOLUB, G.H.},
Some modified matrix eigenvalue problems,
{\it SIAM Rev. 15} (1973), 318--334.
\item
{\tt GOLUB, G.H.} and {\tt WELSCH, J.H.},
Calculation of Gauss quadrature rules,
{\it Math. Comp. 23} (1969), 221--230.
\item
{\tt GRAGG, W.B.} and {\tt HARROD, W.J.},
The numerically stable reconstruction of Jacobi matrices from
spectral data, {\it Numer. Math. 44} (1984), 317--335.
\item
{\tt KAUTSKY, J.} and {\tt GOLUB, G.H.},
On the calculation of Jacobi matrices, {\it Linear Algebra Appl.
52/53} (1983), 439--455.
\item
{\tt LUKE, Y.L.},
{\it Mathematical Functions and Their Approximations}, Academic
Press, New York, 1975.
\item
{\tt REES, C.J.}, Elliptic orthogonal polynomials, {\it Duke Math.
J. 12} (1945), 173--187.
\item
{\tt RUTISHAUSER, H.},
On Jacobi rotation patterns, in {\it Experimental Arithmetic, High
Speed Computing and Mathematics} (N.C. Metropolis et al., eds.),
pp. 219--239. {\it Proc. Symp. Appl. Math. 15}, American Math.
Soc., Providence, R.I., 1963.
\item
{\tt SACK, R.A.} and {\tt DONOVAN, A.F.},
An algorithm for Gaussian quadrature given modified moments, {\it
Numer. Math. 18} (1971/72), 465--478.
\item
{\tt STROUD, A.H.} and {\tt SECREST, D.},
{\it Gaussian Quadrature Formulas},\linebreak Prentice-Hall, Englewood
Cliffs, N.J., 1966.
\item
{\tt UVAROV, V.B.},
Relation between polynomials orthogonal with different weights
(Russian), {\it Dokl. Akad. Nauk SSSR 126} (1959), 33--36.
\item
{\tt UVAROV, V.B.},
The connection between systems of polynomials that are orthogonal
with respect to different distribution functions (Russian),
{\it \v{Z}. Vy\v{c}isl. Mat. i Mat. Fiz. 9} (1969),
1253--1262.
\item
{\tt WHEELER, J.C.},
Modified moments and Gaussian quadrature,
{\it Rocky Mountain J. Math. 4} (1974), 287--296.
\item
{\tt WILKINSON, J.H.},
{\it The Algebraic Eigenvalue Problem}, Clarendon Press, Oxford,
1965.
\end{enumerate}
\end{document}